\documentclass[12pt]{article}

\catcode`\ =\active \newcommand {\`e} \catcode`\ˆ =\active \newcommand ˆ{\`a} \catcode`\ =\active \newcommand {\`u}
\catcode`\Ž =\active \newcommand Ž{\'e} \catcode`\ =\active \newcommand {\c c} \newcommand \parn{\par\noindent} 
\catcode`\ =\active \newcommand {\^e} \catcode`\™ =\active \newcommand ™{\^o} 
\catcode`\" =\active \newcommand "{\^\i } \catcode`\ž =\active \newcommand ž{\^u} 

\catcode`\: =\active \newcommand :{\ifmmode\string:\else\unskip\kern 2pt\string: \ignorespaces\fi } 
\catcode`\; =\active \newcommand ;{\ifmmode\string;\else\unskip\kern 2pt\string; \ignorespaces\fi } 

\font \db = msbm10 at 12 pt
\font \sdb = msbm7 at 8.5 pt     \font \bdb = msbm12 at 14 pt   \font \Bdb = msbm16 at 18 pt  

     \def \bH{\mbox{\bdb H}}  \def\BR{\mbox{\Bdb R}} 

\renewcommand\a{\alpha}    \newcommand\g{\gamma}    \newcommand\D{\Delta}  \newcommand \G{\Gamma}
\newcommand \la{\lambda}  \renewcommand\o{\omega}    \renewcommand\t{\theta}  
\newcommand\e{\varepsilon}  \newcommand\f{\varphi}    \newcommand\rr{\varrho}   

\newcommand\R{\mbox{\db R}}  \renewcommand\H{\mbox{\db H}} \newcommand\E{\mbox{\db E}} \renewcommand\P{\mbox{\db P}} 
\newcommand\N{\mbox{\db N}}    
\newcommand\sR{\mbox{\sdb R}} \newcommand\sH{\mbox{\sdb H}} 
  
\newcommand\sN{\mbox{\sdb N}} 

\renewcommand\S{\mbox{\db S}}      \newcommand\Z{\mbox{\db Z}}  
\newcommand\Q{\mbox{\db Q}}      
  \newcommand\sC{\mbox{\sdb C}}   \newcommand\sZ{\mbox{\sdb Z}} 

    \newcommand\CC{{\cal C}}     
\newcommand\FF{{\cal F}}  \newcommand\GG{{\cal G}}    \newcommand\II{{\mathcal I}}   
\newcommand\LL{{\mathcal L}}     \newcommand\NN{{\cal N}}    \renewcommand\O{{\mathcal O}}

                 \newcommand\xt{\tilde{x}}
\newcommand\yt{\tilde{y}}

\newcommand\p{\partial}        \newcommand\ii{\infty}  \newcommand\sm{\setminus}  
\newcommand\ra{\rightarrow}  \newcommand\lra{\longrightarrow}       \newcommand\Ra{\Rightarrow}  
  \newcommand\LRa{\Leftrightarrow}    \newcommand\LRA{\Longleftrightarrow}
                 \newcommand\sea{\searrow}

\newcommand\ub{\underbar}      \newcommand\ol{\overline} 

\newcommand\bigint{{\ds\int}}   \newcommand\ds{\displaystyle}   \newcommand\ts{\textstyle}       \renewcommand\ss{\scriptstyle}
\def\rt1{\sqrt{-1}\,\,}  
\newcommand\1{^{-1}}  \newcommand\2{^{-2}}  \newcommand\5{{\ts {1\over 2}}}
\renewcommand\parn{\par\noindent} \newcommand\moins{{\ss \sm}}    \newcommand\indf{\leavevmode\indent }     

   \newcommand\tg{{\rm tg}\,}           \renewcommand\th{{\rm th}\,}        
\newcommand\ch{{\rm ch}\,}       \newcommand\sh{{\rm sh}\,}   \newcommand\arctg{{\rm arctg}\,}
 \newcommand\dist{{\rm dist}\,}   \def\Re{{\cal R}\!{\it e}}    \def\Im{{\cal I}\!{\it m}}

\begin{document}

\newtheorem{defi}{Definition}[section]
\newtheorem{theo}{Theorem}[section]
\newtheorem{prop}{Proposition}[section]  \newtheorem{propr}{Property}[section]
\newtheorem{cor}{Corollary}[section]
\newtheorem{lem}{Lemma}[section]
\newtheorem{rem}{Remark}[section]

\newcommand \beq{\begin{equation}} \newcommand \eeq{\end{equation}} 

\newcommand \bthe{\begin{theo}}  \newcommand \ethe{\end{theo}}   
\newcommand \bpro{\begin{prop}}  \newcommand \epro{\end{prop}}   
\newcommand \bcor{\begin{cor}}    \newcommand \ecor{\end{cor}}     
\newcommand \blem{\begin{lem}}   \newcommand \elem{\end{lem}}   
\newcommand \brem{\begin{rem}}  \newcommand \erem{\end{rem}} 
\newcommand \bdefi{\begin{defi}}   \newcommand \edefi{\end{defi}}

\title{\bf Asymptotic homology of the quotient of $PSL_2(\BR)$ by a modular group}

\author{Jacques FRANCHI}
\medskip 
\date { I.R.M.A., Universit\'e Louis Pasteur et CNRS, 7 rue Ren\'e
Descartes, \newlineÊ\newlineÊ 
67084 \  Strasbourg cedex. \quad France } 
\maketitle
\vspace{-5mm} 
\ub{Electronic mail} :\quad  franchi@math.u-strasbg.fr  
\bigskip 

\begin{abstract}
Consider $\,G:= PSL_2(\R)\equiv T^1\H^2$, a modular group $\,\Gamma$, and the 
homogeneous space $\,\Gamma\!\sm\! G \equiv T^1(\Gamma\!\sm\!\H^2)$. Endow $\,G\,$, and then $\,\Gamma\!\sm\! G\,$, with a canonical left-invariant metric, thereby equipping it with a quasi hyperbolic geometry. \   Windings around handles and cusps of $\,\Gamma\!\sm\! G\,$ are calculated by integrals of closed 1-forms of $\,\Gamma\!\sm\! G\,$. The main results express, in both Brownian and geodesic cases, the joint convergence of the law of these integrals, with a stress on the asymptotic independence between slow and fast windings. \     The non-hyperbolicity of $\,\Gamma\!\sm\! G $ is responsible for a difference between the Brownian and geodesic asymptotic behaviours, difference which does not exist at the level of the Riemann surface $\Gamma\!\sm\!\H^2$ (and generally in hyperbolic cases). \     Identification of the cohomology classes of closed 1-forms and with harmonic 1-forms, and equidistribution of large geodesic spheres, are also addressed. 
\end{abstract}

\bigskip 
\noindent 
\ub{Keywords} :\quad  Brownian motion, Geodesics, Geodesic flow, Ergodic measures, \par Asymptotic laws, Modular group, Quasi-hyperbolic manifold, Closed 1-forms. \par \medskip 
\noindent 
\ub{Mathematics Subject Classification 2000} : \quad primary \   58J65 ; \par 
secondary \  60J65, 37D40, 37D30, 37A50, 20H05, 53C22. \par \medskip 
\noindent  
\ub{Abbreviated title} :\quad  Asymptotic homology of a modular quotient. \par 

\eject

\section{Introduction} \label{sec.Introd} \indf 
    Consider $\,G:= PSL_2(\R)\equiv T^1\H^2$, a modular group $\,\Gamma$, and the 
homogeneous space $\Gamma\!\sm\! G \equiv T^1(\Gamma\!\sm\!\H^2)$. Endow $\,G\,$, and then $\,\Gamma\!\sm\! G\,$, with a canonical left-invariant metric, thereby equipping it with a quasi-hyperbolic geometry, which pertains to the 6th 3-dimensional geometrical structure of the eight described by Thurston [T]. The non-hyperbolic manifold $\,\Gamma\!\sm\! G\,$ has finite volume, finite genus, and a finite number of cusps.  \par 
   It is natural in this setting to study the asymptotic behaviour of the Brownian motion and of the geodesic flow (under some Liouville-like measure), by means of their asymptotic homology, calculated by the integrals of the harmonic 1-forms of $\,\Gamma\!\sm\! G$ along their paths. The main results here express, in both Brownian and geodesic case, the joint convergence in law of these integrals, with asymptotic independence of slow and fast windings. \par 
   This same results yield in fact also the asymptotic law of the normalised integrals, along Brownian and geodesic paths, of  any $C^2$ closed 1-form. Indeed, it holds true on $\,\Gamma\!\sm\! G$ that the cohomology classes of closed 1-forms can be identified with harmonic 1-forms. \par 
   The non-hyperbolicity of $\,\Gamma\!\sm\! G\,$ is responsible for a difference between the Brownian and geodesic asymptotic behaviours, difference which does not exist at the level of the Riemann surface $\Gamma\!\sm\!\H^2$. Lifting to the unit tangent bundle has also the advantage to free the harmonic forms of the constraint to have a null sum of their residues at the cusps of $\Gamma\!\sm\!\H^2$. Counter to the hyperbolic setting, the geodesic flow on $\,\Gamma\!\sm\! G\,$ is not ergodic with respect to the normalised Liouville measure of $\,\Gamma\!\sm\! G$, so that the Liouville-like law governing the geodesic, for which ergodicity holds, has to be supported by some leaf. Moreover asymptotic equidistribution of large geodesic spheres holds for such measure. \par \smallskip 

   This article, in which the modular group $\,\Gamma$ is arbitrary, is mainly a generalisation of [F3], which deals with the particular case of the modular group $\,\Gamma$ being the commutator subgroup of $\,PSL_2(\Z)$ ; in which case the quotient  manifold (which is interestingly linked to the trefoil knot) has a unique cusp and a unique handle ; for example, we have now to take into account 
the pullbacks of singular harmonic forms on $\Gamma\!\sm\!\H^2$, which did not exist in [F3]. \  
  However, the identification of the cohomology classes with harmonic spaces $\,H^1$ and the equidistribution of large geodesic spheres, two questions which are addressed here, were not discussed in [F3]. \par  \smallskip 
  
     Brownian and geodesic asymptotic behaviours were already studied in a similar way, but in an hyperbolic setting, in [E-F-LJ1], [E-F-LJ2], [E-LJ], [F2], [G-LJ], [LJ1], [LJ2], [W].  \   As in [F3], hyperbolicity (which is replaced by quasi-hyperbolicity) does not hold in the present setting, nor ergodicity of the Liouville measure (which has to be replaced by Liouville-like measures, supported by leaves), and the asymptotic Brownian and geodesic windings are no longer the sames, though
comparable (the spiral windings of the geodesics of $\,G\,$ about their projections on $\H^2$ is mainly responsible for this feature). Moreover fast and slow windings are addressed jointly. \par\smallskip 

   As in  [LJ1], [G-LJ], [F3], the aim is here to study asymptotic homology, meaning that only closed forms are considered. Consequently no foliated diffusion is needed. Whereas [LJ2],  [E-LJ], [F2], [E-F-LJ1], [E-F-LJ2] dealt with non necessarily everywhere closed 1-forms, so that the showing up of a spectral gap at the level of the stable foliation was needed. 
   
\subsection{Outline of the article} \label{sec.outline} \indf 
   The framework of this article is along the following sections, as follows. 
\smallskip \par \noindent 
{\bf \ref{sec.Iwas}) \  {Iwasawa coordinates and metrics on ${\ts G= PSL_2(\R )}$}} \par 
   Taking advantage of the global Iwasawa coordinates on $\,G\,$, a canonical one-parameter family of left-invariant Riemannian metrics on $\,G\,$ is exhibited, endowing it with a non-hyperbolic, but quasi-hyperbolic geometry (of 3-dimensional tangent bundle). \smallskip \parn  
{\bf \ref{sec.Geom}) \  {Geometry of a modular homogeneous space $\,\Gamma\!\sm\! G$} }  \par 
   A basis of harmonic 1-forms on $\,\Gamma\!\sm\! G$ is described in Theorem \ref{the.cohom}, together with their asymptotics in the cusps. A particular role is played by a form $\,\omega_0\,$, which is not the pullback of a form on $\,\Gamma\!\sm\!\H^2$. Comparing with the two dimensional case of $\Gamma\!\sm\!\H^2$, lifting to the unit tangent bundle $\,\Gamma\!\sm\! G$ has then also the advantage to free the harmonic forms of the constraint to have a null sum of their residues at the cusps. 
\smallskip \parn 
{\bf \ref{sec.closform}) \ {Closed forms and harmonic forms} }  \par 
   The geometries associated with $\,\Gamma\,$ and with a free normal subgroup of finite index $\,\tilde\Gamma\,$ are compared. The identification of closed 1-forms with harmonic 1-forms modulo exact forms is deduced. This has the important consequence that  the asymptotic study of (integrals of) closed 1-forms will reduce to the asymptotic study of harmonic 1-forms. \smallskip \parn   
{\bf \ref{sec.Brown}) \ {Left Brownian motion on ${\ts G= PSL_2(\R )}$} }  \par 
   The natural left Brownian motion is seen to decompose into a planar hyperbolic Brownian motion and a correlated angular Brownian motion. 
\smallskip \parn   
{\bf \ref{sec.exc}) \ {Asymptotic Brownian windings in $\,\Gamma\!\sm\! G$} }  \par 
   The harmonic forms of the basis $(\omega_j , \tilde\omega_l)_{j,l}\,$ exhibited by Theorem \ref{the.cohom} are integrated along the Brownian paths, run during a same time $\,t\,$ going to infinity. This yields on one hand slow martingales $(\tilde M_t^l)$, accounting for the Brownian windings around the handles, and on the other hand fast martingales $(M_t^j)$, accounting for the Brownian windings around the cusps. \  Theorem  \ref{the.CV} gives the joint asymptotic law of all these normalised martingales. Its statement is mainly as follows :  \parn 
Theorem  \ref{the.CV} {\it \quad    ${\ds \Big( {M^j_t\Big/ t}\, ,\, {\tilde M^l_t\Big/\sqrt{t}\,}
\Big)_{j,l}\, } $ converges in law towards ${\ds \bigg( \sum_{\ell=1}^{\nu_\ii(\Gamma)} r^\ell_j\, {\mathcal Q}_\ell\, , {\mathcal N}^l \bigg)_{j,l}\,}$, where all variables $\,{\mathcal Q}_\ell ,{\mathcal N}^l\,$ are independent, each $\,{\mathcal Q}_\ell\,$ is Cauchy with parameter $\,{h_\ell\over 2\,V(\Gamma\sm\sH^2)}\,$, $\,{\mathcal N}^l$ is centred Gaussian with variance $\langle\tilde\omega_l , \tilde\omega_l\rangle $, and the $\,r^\ell_j\,$ are residues in the cusps.  } 
\par \medskip  \parn   
{\bf \ref{sec.Geodesics}) \ {Geodesics of $\,{\ts G= PSL_2(\R )}$} and ergodic measures}  \par 
   The geodesics of $\,G\,$ are described. They project on $\,\H^2$ as quasi-geodesics having constant speed. Thus $\,T^1(\Gamma\!\sm\! G)$ appears as naturally foliated, with on each leaf an ergodic measure, image of the Liouville measure  of $\,\Gamma\!\sm\! G\equiv T^1(\Gamma\!\sm\! \H^2)$, introduced in Definition \ref{def.Liou}. On the contrary, in this non-hyperbolic structure, the geodesic flow is not ergodic with respect to the Liouville measure  on $\,T^1(\Gamma\!\sm\! G)$. 
\smallskip \parn   
{\bf \ref{sec.Agw}) \ {Asymptotic geodesic windings} }  \par 
   The martingales analysed in Theorem  \ref{the.CV} are in this section replaced by the integrals of the same forms $(\omega_j , \tilde\omega_l)_{j,l}\,$, but along geodesic segments (of length $\,t$) instead of Brownian paths. The geodesics are chosen according to the natural ergodic measures introduced in Section \ref{sec.Geodesics}. \  Theorem \ref{the.geo} describes the asymptotic law of the normalised geodesic windings produced in this way. Its statement is mainly as follows : \smallskip \parn 
Theorem \ref{the.geo}  {\it \quad  ${\ds \Big( t^{-1}\!\int_{\g [0,t]}\o_j\, ,\, t^{-1/2}\! \int_{\g [0,t]}\tilde\o_l\,\Big)_{j,l} } $ \    converges in law, under the ergodic measure $\mu^{k}_\e (d\g )$,  to   ${\ds \left(  {(1+a^2) k \over \sqrt{1+a^2k^2}} \, 1_{\{ j=0\}} + 2 \Big({1-k^2\over 1+a^2k^2}\Big)^{1/2} \sum_{\ell=1}^{\nu_\ii(\Gamma)} r^\ell_j\, {\mathcal Q}_\ell \, , \Big({4(1-k^2) \over 1+a^2k^2}\Big)^{1/4} {\mathcal N}^l \right)_{j,l}} , $  
where the limit random variables ${\mathcal Q}_\ell , {\mathcal N}^l$ are as in Theorem  \ref{the.CV}, and $\,a\,$ is the parameter of the metric. } \par\smallskip 
      An interesting feature is the difference between the Brownian and geodesic behaviours, in noteworthy contrast with the hyperbolic case : counter to the Brownian case, the $\,d\t$-part of the form $\,\o_0\,$ is responsible for a non-negligible asymptotic contribution, and the metric parameter $\,a\,$ now appears in the limit law. 
\smallskip \parn   
{\bf \ref{sec.Erlgs}) \ {Equirepartition in $\,\Gamma\!\sm\! G\,$ of large geodesic spheres} }  \par 
    Corollary \ref{cor.Erlgs} asserts that the ergodic measures $\,\mu^{k}_\e\,$ of Section \ref{sec.Geodesics} and Theorem \ref{the.geo}, are weak limit of the uniform law on large geodesic quasi-spheres. This is easily deduced from the following. \parn 
Theorem \ref{the.Erlgs} \  {\it  The normalized Liouville measure $\,\mu^\Gamma$ on $\,T^1(\Gamma\!\sm\!\H^2)\equiv \Gamma\!\sm\! G\,$ is the weak limit as $\,R\to\ii\,$ of the uniform law on the geodesic sphere  $\,\Gamma g\,PSO(2) \Theta_{R}$ of $\,\Gamma\!\sm\! G\,$ having radius $\,R\,$ and fixed center $\,\Gamma g\in \Gamma\!\sm\! G$ : for any  compactly supported continuous function $\,f\,$ on $\,\Gamma\!\sm\! G\,$, denoting by $\,d\rr\,$ the uniform law on $PSO(2)$, we have 
$$  \int f\,d\mu^\Gamma = \lim_{R\to\ii} \int_{PSO(2)} f(\Gamma g\,\rr \,\Theta_{R})\, d\rr \, . $$ 
} 
 The equidistribution theorem \ref{the.Erlgs} and its proof (based on the mixing theorem) were  already given by Eskin and McMullen in [E-MM].  
\smallskip \parn   
{\bf \ref{sec.SyntPr}) \ {Synthetic proof of Theorem \ref{the.CV}} }  \par 
    This proof is an adaptation of an analogous proof in ([F3], Section 10). Thus some details are here somewhat eluded, for which we refer to [F3]. However all ingredients are given, with a stress on differences with the particular case addressed in [F3], and the most involved arguments are detailed to a certain extent.  \par     
    The main difficulty of the whole proof, widely responsible for its length, is to establish the asymptotic independence between slow windings (about the handles) and singular windings (about the cusps). This demands in particular to get good approximation of the contribution of both type, and then to analyse carefully the successive excursions of Brownian motion in the core and in the cusps of the quotient hyperbolic surface. \par\smallskip \parn 
{\bf \ref{sec.proofg}) \ {Proof of Theorem \ref{the.geo}} }   \par 
   The strategy for this proof is mainly to replace the geodesic paths by the Brownian paths, as in hyperbolic case ([LJ1], [LJ2], [E-LJ], [F2], [E-F-LJ1], [E-F-LJ2]), and as in [F3], in order to reduce Theorem \ref{the.geo} to Theorem \ref{the.CV}. 
Here again, the analogous proof in ([F3], Section 14) is adapted. 

\section{Iwasawa coordinates and metrics on ${\ts G= PSL_2(\BR )}$} \label{sec.Iwas} 
This section is mainly taken from [F3].  \par   
   Consider the group $G := PSL_2(\R )$, which is classically parametrized by the Iwasawa coordinates \   $ (z=x+\rt1 y\, ,\, \t )\in\H^2\times (\R /2\pi\Z )\,$ ($\H^2$ denotes as usual the hyperbolic plane, identified with the PoincarŽ half-plane), in the following way : each $g\in G$ writes uniquely 
$$ g=g(z,\t ) := \pm\, n(x) a(y) k(\t )\,, $$
where $\; n(x)\, ,\, a(y)\, ,\, k(\t )\;$ are the one-parameter subgroups defined by : 
\beq \label{f.semigr} n(x) := \pmatrix{ 1 & x \cr 0 & 1 \cr} , \; a(y) := \pmatrix{ \sqrt{y} & 0 \cr 0 & 1/\sqrt{y} \cr} , \;  k(\t ) := \pmatrix{ \cos (\t /2) & \sin (\t /2) \cr -\sin (\t /2) & \cos (\t /2) \cr}  , \eeq 
and generated respectively by the following elements of the Lie algebra $\, s\ell_2(\R )$ : 
$$ \label{f.sl2R} \nu := \pmatrix{ 0 & 1 \cr 0 & 0 \cr} , \; \a := \pmatrix{ 1/2 & 0 \cr 0 & -1/2 \cr} ,
\;  \kappa := \pmatrix{ 0 & 1/2 \cr -1/2 & 0 \cr} . $$ 
Note that \qquad $ g=g(z,\t ) \LRA \Big[ g(\rt1 ) = z \;\; \hbox{ and } \;\; g'(\rt1 ) =\, y\,
e^{\rt1\t}\,\Big]\, . $ \par\medskip 

   Set also \ $\;\la := \nu -\kappa = \pmatrix{ 0 & 1/2 \cr 1/2 & 0 \cr}\,  $, \ which is natural,
since $\,\a ,\la\,$ are symmetrical while $\,\kappa\,$ is skew-symmetrical, and since in the basis $\; (\alpha ,\la ,\kappa )\;$ of $\, s\ell_2(\R )\,$ the Killing form is diagonal : it has matrix 
$\;\pmatrix{ -2 & 0 & 0 \cr 0 & -2 & 0\cr 0 & 0 & 2\cr}\, $.
\par\smallskip 

   For this reason, we take on $\, s\ell_2(\R )\,$ the inner product such that the basis $\; (\a ,\la ,
a\kappa )\;$ is orthonormal, for some arbitrary parameter $a\in\R^*$. And since we want to work on an homogeneous space $\,\Gamma\!\sm\! G\,$, the Riemannian metric to be considered on $G$ must be a least $\Gamma$-left-invariant, and then a natural choice for the Riemannian metric on $G$ is the left-invariant metric, say $\, ((g^a_{ij}))\,$, generated by the above inner product on $\, s\ell_2(\R )\,$. \par  

   The simple lemma below shows that this choice of metric(s) is geometrically canonical (up to a trivial multiplicative constant), $G$ being seen as $T^1\H^2$. This equips $\,G\equiv T^1\H^2$, and its homogeneous spaces as well, with the 6th of the eight 3-dimensional geometries described by Thurston ([T]), and actually with a quasi-hyperbolic but not hyperbolic structure.  \par 

   Let us denote by $\;\LL_\nu\; ,\; \LL_\a\; ,\; \LL_\kappa\; ,\; \LL_\la\;$ the left-invariant vector fields
on $G$ generated respectively by $ \nu \, ,\, \a\, ,\, \kappa\, ,\, \la\,$. \ A standard computation shows that 
$$ \LL_\la = y\sin\t{\p\over\p y} + y\cos\t{\p\over\p x} -\cos\t {\p\over\p \t} \; ,\; 
\LL_\a = y\cos\t{\p\over\p y} - y\sin\t{\p\over\p x} + \sin\t{\p\over\p \t} \; ,\;
\LL_\kappa = {\p\over\p \t} \; . $$ 

\begin{lem} \label{lem.nat}\quad The Riemannian metrics $\, ((g^a_{ij}))\,$ defined above are, up to a multiplicative constant, the only ones on $G$ which are left-invariant and also invariant with respect to the action of the (Cartan compact subgroup) circle $\,\exp(\R\kappa) = \{ k(\theta )\}\,$. They are given in Iwasawa coordinates $\, (y,x,\t )\,$ by : 
\beq \label{f.gaij} ((g_{i j}^a)) := \pmatrix{ y^{-2} & 0 & 0 \cr 0 & (1+a^{-2})y^{-2} & a^{-2}y\1\cr 0 & a^{-2}y\1 & a^{-2}\cr} . \eeq 
\end{lem} 
\ub{Proof} \quad The left-invariant metrics on $G$ are those which are given by a constant matrix 
$\, ((a_{ij}))\,$ in the basis ${\ds\LL :=  (\LL_\a ,\LL_\la ,\LL_\kappa )}\,$. Set 
$\,{\ds \II := \Big({\p\over\p y} ,{\p\over\p x},{\p\over\p\t}\Big)}\,$. We have $\, \II = \LL A\,$, with 
${\ds A:= \pmatrix{ y^{-1}\cos\t & -y^{-1}\sin\t & 0\cr y^{-1}\sin\t & y^{-1}\cos\t & 0\cr 0 & y\1 & 1\cr}
}$,   so that the left-invariant metrics are given in the basis $\,\II\,$ by $\,^t\! A((a_{ij})) A\,$.  Among
them, the ones we want have to satisfy  $\, {\ds{\p\over\p\t}\, ^t\! A((a_{ij})) A =0}\,$. \parn  
A direct computation shows that this is equivalent to 
$\, ((a_{ij})) = c^2 \pmatrix{ 1 & 0 & 0 \cr 0 & 1 & 0\cr 0 &0 & a^{-2}\cr}$, and then to $((g_{i j}^a))$ being as in the statement. $\;\diamond $ 
\par\medskip 

   Note that with these metrics any holomorphic form $\;f(z) dz\;$ is coclosed, and then harmonic. \par \smallskip 

   The left Laplacian on $G$ corresponding to the basis $\; (\a ,\la , a\kappa )\;$ is the Beltrami Laplacian associated with the  metric $ ((g_{i j}^a)) $, and is given by 
\beq \label{f.Lapl}  \D^a := \LL_\la^2 + \LL_\a^2 + a^2\LL_\kappa^2 = y^2\, \Big({\p^2\over\p y^2}+
{\p^2\over\p x^2}\Big) - 2y\, {\p^2\over\p \t\p x} + (1+a^2)\, {\p^2\over\p \t^2} \; . \eeq  

   Note that $\, \LL_\la\,$ and $\,\LL_\a\,$ generate the canonical horizontal left-invariant vector fields lifted from $\H^2$ to $G$, $\H^2$ being endowed with its Levi-Civita connexion, so that $\;\D^0\,$ is the Bochner horizontal left Laplacian, and $\, \D^a = \D^0 + a^2 {\p^2\over\p \t^2} \;$. \par 

   The measure ${\ds\; \mu (dg) := { dx\, dy\, d\t\over 4\pi^2\, y^2}}\;$ is bi-invariant, hence this is both the Haar measure of $\,G\,$ and the Liouville measure of $\,T^1\H^2$. \par\medskip
   
   Recall that an isometry $\,\gamma\in G\,$ is called respectively elliptic, parabolic, or loxodromic, according as it fixes a point in $\,\H^2$, no point in $\,\H^2$ and a unique point in $\,\partial\H^2=\R\cup\{\ii\}$, or no point in $\,\H^2$ and two points in $\,\partial\H^2=\R\cup\{\ii\}$, respectively. Any isometry of $\,\H^2$ is either elliptic, or parabolic, or loxodromic.  We shall use the following easy lemma. \par 
\blem \label{lem.exp} \quad Any parabolic or loxodromic isometry $\,\gamma\in G\,$ can be written
$\,\gamma = \pm\exp(\sigma)$, for a unique $\,\sigma \in sl_2(\R)$. 
\elem
\ub{Proof} \quad  Consider first a loxodromic $\,\gamma\in G\,$, and an isometry $\,g\in G\,$ mapping two fixed boundary points of $\,\gamma\,$ to $\{ 0,\ii\}$, so that we have for some $\,t\in\R^*$ :  
$$ \gamma = \pm\, g\1\pmatrix{e^t&0\cr 0&e^{-t}} g = \pm \exp\left[ g\1\pmatrix{t&0\cr 0&{-t}} g\right] 
. $$ 
And if ${\ds\,\sigma = \pmatrix{u&v+w\cr v-w&{-u}} \in sl_2(\R)}$, then \ ${\ds \sigma^2= (u^2+v^2-w^2) \pmatrix{1&0\cr 0&1}}$, and then  
$$ \exp(\sigma) = (\ch\rr)\, \pmatrix{1&0\cr 0&1} + {\sh\rr\over \rr}\, \, \sigma \quad \hbox{or} \quad 
 \exp(\sigma) = (\cos\rr)\, \pmatrix{1&0\cr 0&1} + {\sin\rr\over \rr}\, \, \sigma \, , $$ 
according as \  $ (u^2+v^2-w^2) =: \pm\rr^2\,$ is non-negative or negative. \par 
   Hence (setting $\,\sigma := g\sigma'g\1$) : \  ${\ds \gamma = \pm\exp(\sigma') \LRa \pmatrix{e^t&0\cr 0&e^{-t}} = \pm\exp(\sigma)}\,$, which in the first case implies at once $\,v=w=0\,$,  whence $\,u=t\,$, and in the second case : $\,\sin\rr =0\,$, whence $\,e^t=e^{-t}$, an impossibility, establishing the unicity of $\,\sigma'$.  \par 
   Consider then a parabolic $\,\gamma\in G\,$, and an isometry $\,g\in G\,$ mapping its fixed boundary point to $\,\ii\,$, so that we have for some $\,x\in\R^*$ :  
$$ \gamma = \pm\, g\1\pmatrix{1&x\cr 0&1} g = \pm \exp\left[ g\1\pmatrix{0&x\cr 0&0} g\right] . $$ 
And \  ${\ds \gamma = \pm\exp(\sigma') \LRa \pmatrix{1&x\cr 0&1} = \pm\exp(\sigma)}\,$, which in the first case (for $\,\exp(\sigma)$) implies at once $\,u=v-w=0\,$,  whence $\,v=x/2\,$, and in the second case : $\,\sin\rr =0\,$, whence an impossibility, establishing again the unicity of $\,\sigma'$.  $\;\diamond$    

\section{Geometry of a modular homogeneous space $\,\Gamma\!\sm\! G$} \label{sec.Geom} \indf 
    Consider the group $G := PSL_2(\R )$, its full modular subgroup $\Gamma(1) := PSL_2(\Z)$, and another modular subgroup $\,\Gamma$, that is, a subgroup of $\,\Gamma(1)$ having finite index $\,[\Gamma(1):\Gamma]$. \par 
       As usual, let us identify $\,G\,$ with the unit tangent bundle $T^1\H^2\equiv \H^2\times\S^1$ of the hyperbolic plane $\H^2$, and also with the group of M\"obius isometries
(homographies $\,z\mapsto {az+b\over cz+d}$ with $\,ad-bc =1$) of $\H^2$, that is the group of direct isometries of $\H^2$. \par 
   The elements $\, u :=(z\mapsto -1/z)\, ,\, v:=(z\mapsto (z-1)/z)\,$ generate the group $\Gamma(1)$, which admits the presentation $\{ u,v\, |\, u^2=v^3=1\}$. $\Gamma(1)$ is of course also generated by $\{ u,vu=(z\mapsto z+1)\}$. \par \smallskip 

    Note that $\,[\Gamma,\Gamma]\,$ is a free group, as a subgroup of $\,D\Gamma(1) := [\Gamma(1),\Gamma(1)]$, which is the free group generated by ${\ds \pm \pmatrix{2&1\cr 1&1}}\,$ and ${\ds \pm \pmatrix{1&1\cr 1&2}}\,$.  \   Generally $\,[\Gamma : [\Gamma,\Gamma]] $ has not to be finite, as shows the counterexample  $\,D\Gamma(1)$, for which $\,D\Gamma(1) / [D\Gamma(1),D\Gamma(1)]\equiv \Z^2$. \  But \parn
since $\,D\Gamma(1)$ is a normal subgroup of $\,\Gamma(1)$ such that $\,\Gamma(1)/D\Gamma(1) \equiv \Z/6\Z\,$, then 
$$ \tilde\Gamma := \Gamma\cap D\Gamma(1) $$ 
is a free and normal subgroup of  $\,\Gamma$ such that $\,\Gamma / \tilde\Gamma \equiv \Gamma\cdot D\Gamma(1)\Big/ D\Gamma(1)$ is a subgroup of $\,\Gamma(1)/D\Gamma(1)$, so that \  $\Gamma / \tilde\Gamma$ is isomorphic to a subgroup of $\,\Z/6\Z\,$.  \par \smallskip 

\medskip

   We are interested in the modular homogeneous space $\,\Gamma\!\sm\! G\,$. The identification of  $\,G\,$ with the unit tangent bundle $T^1\H^2$ allows to identify similarly this modular homogeneous space with the unit tangent bundle of the corresponding Riemann surface $\Gamma\!\sm\! \H^2$ :  
$$\Gamma\!\sm\! G \,\equiv\, \Gamma\!\sm T^1\H^2 \,\equiv\, T^1(\Gamma\!\sm\! \H^2) . $$ 

   The projection of the Liouville measure ${\ds\; \mu (dg) = { dx\, dy\, d\t\over 4\pi^2\, y^2}}\;$ 
(which is also a right and left Haar measure on $\,G$) onto $\,\Gamma\!\sm\! G\,$ is  proportional to the volume measure $\,V\,$ of $\,\Gamma\!\sm\! G\,$. \  
By the choice of the metric $((g^a_{ij}))$, the volume of $\,\Gamma\!\sm\! G\,$ is clearly \  $V(\Gamma\!\sm\! G) = {2\pi\over |a|}\times {\rm covol}(\Gamma)$, 
where $\,{\rm covol}(\Gamma) = V(\Gamma\!\sm\!\H^2)$ denotes the finite hyperbolic volume of $\,\Gamma\!\sm\! \H^2$.    \par
   Let \  $\mu^\Gamma := {2\pi\over  {\rm covol}(\Gamma)}\, \mu_{\Big| \Gamma\sm\! G}\,$ denote the normalized projection of $\,\mu\,$ on $\,\Gamma\!\sm\! G\,$, identified with a law on left $\,\Gamma$-invariant functions  on $\,G$.   
   
\subsection{From $\Gamma\!\sm\! \bH^2$ to $\,\Gamma\!\sm\! G $} \label{sec.fromto} \indf  
   The following lemma ensures that the lift to the unit tangent bundle increases the first Betti number of  the Riemann surface $\Gamma\!\sm\! \H^2$ by exactly one. I thank T. Delzant for having explained to me why, so that I owe to him this lemma. 
\blem \label{lem.homol} \quad The modular homogeneous space $\,\Gamma\!\sm\! G \equiv\, T^1(\Gamma\!\sm\! \H^2)$ is diffeomorphic to $(\Gamma\!\sm\! \H^2)\times\S^1$. Consequently, we have the following simple relation between the first Betti numbers of $\,\Gamma\!\sm\! G$ and of $\,\Gamma\!\sm\! \H^2$ : \quad 
$ {\rm dim}\, [H^1(\Gamma\!\sm\! G)] = 1+  {\rm dim}\, [H^1(\Gamma\!\sm\!\H^2)] . $ 
\elem 
\ub{Proof} \quad  As a cover of $\,\Gamma(1)\!\sm\! \H^2$, the Riemann surface $\,\Gamma\!\sm\! \H^2$ is orientable and non-compact. As is known for any orientable non-compact smooth manifold, it carries a smooth non-vanishing vector field, hence a smooth cross section $\Big(x\mapsto (x,\vec{v}_x)\in T^1_x(\Gamma\!\sm\! \H^2)\Big)$ of $\,T^1(\Gamma\!\sm\! \H^2)$. Denoting by $\,\alpha =\alpha_x(\vec{v})$ the angle $\widehat{(\vec{v}_x,\vec{v})}$ in the oriented plane $\,T^1_x(\Gamma\!\sm\! \H^2)$, we get the diffeomorphism : $\,(x,\vec{v})\mapsto (x,\alpha)$ from  $\,T^1(\Gamma\!\sm\! \H^2)$ onto $(\Gamma\!\sm\! \H^2)\times\S^1$. $\;\diamond$  
\medskip 

   Furthermore, under the canonical projection $\,\pi : \Gamma\!\sm\! G \ra \Gamma\!\sm\! \H^2\equiv \Gamma\!\sm\! G/\exp(\R\kappa)$, the harmonic space $\,H^1(\Gamma\!\sm\! \H^2)$ (that is, the space of real harmonic forms on $\,\Gamma\!\sm\! \H^2$) is pulled back to the subspace $\,\pi_*[H^1(\Gamma\!\sm\! \H^2)]$ of the harmonic space $\,H^1(\Gamma\!\sm\! G)$, which is isomorphic to $\,H^1(\Gamma\!\sm\! \H^2)$. \par 
   Thus, to describe the harmonic 1-forms of the modular homogeneous space $\,\Gamma\!\sm\! G\,$, once the harmonic 1-forms  of the Riemann surface $\,\Gamma\!\sm\! \H^2$ are known, by the above lemma \ref{lem.homol} it is sufficient to produce a harmonic 1-form $\,\omega_0\in H^1(\Gamma\!\sm\! G)$, such that $\,\omega_0\notin  \pi_*[H^1(\Gamma\!\sm\!\H^2)]$. \par 
   Now, such harmonic 1-form $\,\omega_0\,$ was computed in ([F3], Section 6), as the restriction to $\,\Gamma\!\sm\! G\,$ of a harmonic 1-form on $\,\Gamma(1)\!\sm\! G$ :
\beq \label{f.omega0} \omega_0 := d\t + 4\, \Im (\eta'(z)/\eta(z))\, dx + 4\, \Re (\eta'(z)/\eta(z))\, dy = d\Big( \t + 4\, \arg (\eta (z))\Big)\, , \eeq 
where $\eta$ denotes the Dedekind function, defined on $\H^2$ (seen as the Poincar\'e half-plane) by : 
$$ \eta (z) := e^{\rt1 \pi\, z/12}\times \prod_{n\in\sN^*}(1-e^{\rt1 2\pi\, n\, z}) \; . $$ 

\subsection{Geometry of the Riemann surface $\Gamma\!\sm\!\bH^2$} \label{sec.cohomG} \indf  
   Let us describe now the harmonic space $\,H^1(\Gamma\!\sm\! \H^2)$ of the Riemann surface $\,\Gamma\!\sm\! \H^2$. Denote by $\,{\bf g}(\Gamma)$ the genus of $\,\Gamma\!\sm\! \H^2$, and by $\,\nu_\ii(\Gamma) := {\rm Card}(\Gamma\!\sm\! \ol{\Q})$ the number of its cusps, that is, the number of $\Gamma$-inequivalent parabolic points of $\,\Gamma$. Note that clearly  $\,\nu_\ii(\Gamma)\ge 1\,$. \par    
   Let us denote these cusps by $\,\CC_1,..,\CC_{\nu_\ii(\Gamma)}\,$, choosing $\,\CC_{\nu_\ii(\Gamma)}=\Gamma \ii\,$ to be the one cusp associated to the particular parabolic point $\,\ii\,$. \par
   The Riemann surface $\,\Gamma\!\sm\! \H^2$ decomposes into the disjoint union of a compact core, which is a compact surface having genus $\,\Gamma\!\sm\! \H^2$ and a boundary made of $\,\nu_\ii(\Gamma)$ pairwise disjoint circles $\,\S^1$, and of $\,\nu_\ii(\Gamma)$ pairwise disjoint ends, each being diffeomorphic to $\,\S^1\times\R_+^*$ and associated with one of the cusps $\,\CC_\ell\,$, which we call ``solid cusp'' and denote also by $\,\CC_\ell\,$.  \par\medskip 
   
   Recall now that the harmonic space $\,H^1(\Gamma\!\sm\! \H^2)$ is the dual of the first real singular homology space $\,H_1(\Gamma\!\sm\! \H^2)$ (see for example ([D], 24.33.2)), so that the above decomposition implies the formula  :
\beq \label{f.H1G} {\rm dim}\, [H^1(\Gamma\!\sm\!\H^2)] = 2\, {\bf g}(\Gamma) + \nu_\ii(\Gamma) - 1\, . \eeq
   
   On the other hand, an automorphic form $\,f\,$ of weight 2 with respect to $\,\Gamma\,$ induces the holomorphic differential $\,f(z)\,dz\,$ on $\,\Gamma\!\sm\! \H^2$. More precisely, let us denote as Miyake ([M]) by $\,\GG_2(\Gamma)$ the complex vector space of those automorphic forms $\,f\,$ (of weight 2) which are holomorphic on $\,\H^2$ and at the cusps of $\,\Gamma$, and by $\,{\cal S}_2(\Gamma)$ the subspace of so-called ``cusp forms'', that is of forms in $\,\GG_2(\Gamma)$ which vanish at the cusps of $\,\Gamma$. \par 
   The so-called Petersson inner product (see for example ([M], Section ¤2.1)) is defined for $\,(f_1,f_2)\in {\cal S}_2(\Gamma)\times \GG_2(\Gamma)$, by : 
\beq \label{f.Peters} \langle f_1(z)dz , f_2(z) dz\rangle = \langle f_1 , f_2\rangle := V(\Gamma\!\sm\! \H^2)\1 \int_{\Gamma\sm\sH^2} f_1(z)\,\ol{f_2(z)}\,dz \, .  \eeq 
Note that $\, y\,|f(z)|\,$ is the natural norm of $\,f(z)dz\in T^*_z(\Gamma\!\sm\! \H^2)$, induced by the volume measure of $\,\Gamma\!\sm\! \H^2$ ; so that the differential $\,\,{\ds |f(z)|^2\,dz = \| f(z) dz\|^2\,{dx dy\over y^2}\,}\,$, integrated over $\,\Gamma\!\sm\! \H^2$, indeed computes precisely the global norm of $\,f(z)dz\in T^*(\Gamma\!\sm\! \H^2)$. \par 
   As in [M] again, let $\,\NN_2(\Gamma)$ denote the orthogonal complement of $\,{\cal S}_2(\Gamma)$  in $\,\GG_2(\Gamma)$, with respect to the Petersson inner product. Then ([M], Theorem 2.5.2) states the following : 
\beq \label{f.GN2} {\rm dim}_{\sC} [ {\cal S}_2(\Gamma) ] = {\bf g}(\Gamma)\quad  \hbox{ and } \quad {\rm dim}_{\sC} [ \NN_2(\Gamma) ] = \nu_\ii(\Gamma) -1\, .  \eeq 
  
   Note that the $-1$ in the second formula is natural, since the sum of the residues of a harmonic differential form has to be zero. \par 
   
   As a consequence, comparing (\ref{f.H1G}) and (\ref{f.GN2}), we see that the regular part of the harmonic space $\,H^1(\Gamma\!\sm\! \H^2)$, that is the part due to the handles, admits a basis made of $\,2\,{\bf g}(\Gamma)$ real harmonic forms $\,\Re [f(z)\, dz]$, with $\,f\in {\cal S}_2(\Gamma)$ a cusp form ; and that the singular part of the harmonic space $\,H^1(\Gamma\!\sm\! \H^2)$, that is the part due to the cusps and having residues at the cusps which are not all null, admits a basis made of $\,(\nu_\ii(\Gamma) -1)$ real harmonic forms $\,\Re [f(z)\, dz]\,$,  with $\,f\in \NN_2(\Gamma)$. \par 
Notice that, specifying that the singular harmonic forms must be real and orthogonal to the regular harmonic forms (and then in particular have non all vanishing residues), we get indeed $\,(\nu_\ii(\Gamma) -1)$ independent such differential forms, and not $\,2 (\nu_\ii(\Gamma) -1)$. 
\par \medskip 

\subsection{Geometry of the modular space $\Gamma\!\sm\! G$} \label{sec.geomG} \indf  
   Note that $\,y\,$ is naturally the height in the solid cusp $\,\CC_{\nu_\ii(\Gamma)}$, or in the corresponding end $\,T^1\CC_{\nu_\ii(\Gamma)}\,$ of $\,\Gamma\!\sm\! G\,$ as well. Similarly, for any solid cusp $\,\CC_\ell\,$, by using a M\"obius isometry mapping it on $\,\CC_{\nu_\ii(\Gamma)}$, we get a $\,\Gamma$-invariant function $\,\yt_\ell\,$ on $\,G\,$ or on $\H^2\,$, which yields a canonical height in the solid cusp $\,\CC_\ell\,$ (or in $\,T^1\CC_{\ell}\,$ as well). Similarly yet, we have canonical coordinates $(\yt_\ell,\xt_\ell,\t)$ on the end $\,T^1\CC_{\ell}\,$, with $\,|d\xt_\ell|=|d\yt_\ell|=\yt_\ell\,$ and $\,|d\t|\equiv 1\,$. \par \smallskip    

   Let us sum up the above, ([F3], Theorem 1)  and ([M], Corollary 2.1.6) in the following. 
\bthe \label{the.cohom} \quad  The real harmonic space $\,H^1(\Gamma\!\sm\! G)$ of the modular homogeneous space $\,\Gamma\!\sm\! G\,$ admits a basis $\Big(\omega_0\,,\,\omega_1, ..,\omega_{(\nu_\ii(\Gamma) -1)}\,,\, \tilde\omega_1, ..,\tilde\omega_{2{\bf g}(\Gamma)}\Big)$, where $\,\omega_0\,$ is given by (\ref{f.omega0}), each $\,\omega_j\,$ (for $1\le j < \nu_\ii(\Gamma)$) equals $\,\Re [f_j(z)\, dz]$ for some automorphic form $\,f_j\in \NN_2(\Gamma)$, and the $\,\tilde\omega_j\,$ are  pairwise orthogonal and can be writen $\,\Re [\tilde f_j(z)\, dz]$ for some cusp form $\,\tilde f_j\,$. \par
   Moreover, the $\,\tilde\omega_j$ are bounded, and we have the following behaviours near the cusps : 
$$ \omega_0 = {\pi\over 3}\, dx + d\t + \O (ye^{-2\pi\, y}) \quad \hbox{ near the cusp } \; \CC_{\nu_\ii(\Gamma)}, \hbox{ id est for }\; y\to\ii\, ; $$
$$ \omega_j = r_j^\ell\, d\xt_\ell + \O (\yt_\ell\,e^{-\pi\,\yt_\ell/h_\ell}) \quad \hbox{ near the cusp } \; \CC_{\ell}\,, \hbox{ id est for }\; \yt_k \to\ii\, ;  $$ 
$r_j^\ell\,$ denoting the residue of the harmonic form $\omega_j $ at the cusp $\CC_\ell$ (so that, in particular, we have for $\,\omega_0$ : $r^{\ell}_0 = {\pi\over 3}\,1_{\{ \ell =\nu_\ii(\Gamma)\}}$), and $\,h_\ell >0\,$ denoting the width of the solid cusp $\,\CC_\ell\,$, determined by : $\,\pm\, n(h_\ell)$ (recall Formula (\ref{f.semigr})) is conjugate in $\,G\,$ to a generator of the parabolic subgroup of $\,\Gamma$ associated with the cusp $\,\CC_\ell\,$.  \  We have \  $\sum_{\ell =1}^{\nu_\ii(\Gamma)}\limits r_j^\ell=0\,$, for $\,1\le j < \nu_\ii(\Gamma)$. 
\ethe 

   Note in particular that, while the regular harmonic forms $\,\tilde\omega_l\,$ belong to $\,L^2(\Gamma\!\sm\! G)$, on the contrary the singular harmonic forms $\,\omega_j\,$ do not belong to $\,L^1(\Gamma\!\sm\! G)$. This is not surprising, since the $\,\tilde\omega_l$ calculate slow windings about the handles of the manifold $\Gamma\!\sm\! G$, whereas the $\,\omega_j$ calculate fast windings about the cusps of the manifold $\Gamma\!\sm\! G$. \par \medskip 

   Note also that lifting to the unit tangent bundle has also the advantage to free the harmonic forms of the constraint to have a null sum of their residues at the cusps of $\Gamma\!\sm\!\H^2$ (which are also the cusps of $\Gamma\!\sm\! G$) : \  the constraint $\sum_{\ell =1}^{\nu_\ii(\Gamma)}\limits r_j^\ell=0\,$, which holds for harmonic (and closed, by Proposition \ref{pro.closedform} below) forms on $\Gamma\!\sm\!\H^2$, does not hold any longer at the level of $\Gamma\!\sm\! G$, since it breaks down in particular for the harmonic form $\,\omega_0\,$.   \par \medskip   
   
   Furthermore, the genus and the volume of $\,\Gamma\!\sm\! \H^2$ can be expressed in terms of two more parameters, the numbers $\,\nu_2(\Gamma)$ and $\,\nu_3(\Gamma)$ of $\,\Gamma$-inequivalent elliptic points of $\,\Gamma$, of order 2 and 3 respectively. \par 
    The genus of $\,\Gamma\!\sm\! \H^2$ is given by the formula ([M], Theorem 4.2.11) : 
\beq \label{f.genus} {\bf g}(\Gamma) = 1+ {\ts{1\over 12}} \,[\Gamma(1) : \Gamma] - {\ts{1\over 4}}\, \nu_2(\Gamma) -{\ts{1\over 3}}\,  \nu_3(\Gamma) - {\ts{1\over 2}}\, \nu_\ii(\Gamma)\, .   \eeq 

    The volume of $\,\Gamma\!\sm\! \H^2$ is given by the formula ([M], Theorem 2.4.3) : 
\beq \label{f.Vol} V(\Gamma\!\sm\! \H^2) = 2\pi\times [ 2 \,{\bf g}(\Gamma) - 2+ \nu_\ii(\Gamma) + \5\, \nu_2(\Gamma) + {\ts{2\over 3}}\, \nu_3(\Gamma) ]\, .   \eeq 

   In the particular case of principal congruence groups $\,\Gamma(N)$, very explicit formulae for $[\Gamma(1) : \Gamma(N)]$ and for $\,\nu_j(\Gamma(N))$ are known (see ([M], Section ¤4.2)), and there exists a precise description of the space $\,\NN_2(\Gamma)$ in terms of analytic continuations of Eisenstein series (see ([M], Section ¤7.2)). \par 
   
\section{Closed forms and harmonic forms}  \label{sec.closform}  \indf 
   Recall that $\,\Gamma\,$ was neither supposed to be a congruence subgroup, nor to be normal. But recall from the beginning of Section \ref{sec.Geom} that \  $\tilde\Gamma := \Gamma\cap D\Gamma(1)$ is a free and normal  subgroup of $\Gamma$, such that $\,\Gamma/\tilde\Gamma$ is isomorphic to a subgroup of $\,\Z/6\Z$. \   Beginning by a comparison between the geometries associated with $\,\Gamma\,$ and $\,\tilde\Gamma\,$, we shall in this section deduce that their cohomology spaces (of 1-forms) identify with their harmonic spaces $\,H^1$, so that the asymptotic study of closed 1-forms will reduce to the asymptotic study of harmonic 1-forms. \par    
   By Formulas (\ref{f.H1G}) and (\ref{f.genus}) we have : 
$$ d := {\rm dim}\, [H^1(\Gamma\!\sm\!\H^2)]  = 2g(\G) + \nu_\ii(\G) -1 \quad \hbox{and} \quad \tilde d := {\rm dim}\, [H^1(\tilde\Gamma\!\sm\!\H^2)] = 2g(\tilde\Gamma) + \nu_\ii(\tilde\Gamma) -1\, , $$ 
$$ {\bf g}(\Gamma) = 1+ {\ts{1\over 12}} \,[\Gamma(1) : \Gamma] - {\ts{1\over 4}}\, \nu_2(\Gamma) -{\ts{1\over 3}}\,  \nu_3(\Gamma) - {\ts{1\over 2}}\, \nu_\ii(\Gamma)\, ; $$ 
and then 
$$ {\bf g}(\tilde\Gamma) = 1+ {\ts{1\over 12}} \,[\Gamma(1) : \tilde\Gamma] - {\ts{1\over 4}}\, \nu_2(\tilde\Gamma) -{\ts{1\over 3}}\,  \nu_3(\tilde\Gamma) - {\ts{1\over 2}}\, \nu_\ii(\tilde\Gamma) $$ 
$$ \ge  1+ \Big( {\ts{1\over 12}} \,[\Gamma(1) : \Gamma]   - {\ts{1\over 2}}\, \nu_\ii(\Gamma) \Big) \times [\Gamma : \tilde\Gamma]  \;\ge\; 1+ [\Gamma : \tilde\Gamma] ({\bf g}(\Gamma)-1)\, , $$ 
since on one hand the free group $\,\tilde\Gamma$ cannot have any elliptic point, and on the other hand, recalling the definition of $\nu_\ii(\Gamma)$ (in Section \ref{sec.cohomG}), we must have : 
$$ \nu_\ii(\tilde\Gamma)= {\rm Card}(\tilde\Gamma\!\sm\! \ol{\Q}) \in[  \nu_\ii(\Gamma) ,  [\Gamma : \tilde\Gamma] \times \nu_\ii(\Gamma)] \, . $$  
From the above we deduce at once :
$$ d = {\rm dim}\, [H^1(\Gamma\!\sm\!\H^2)] = 1+ [\Gamma(1) : \Gamma]/6 - {\ts{1\over 2}}\, \nu_2(\Gamma) -{\ts{2\over 3}}\,  \nu_3(\Gamma) \le 1+ [\Gamma(1) : \Gamma]/6 \, , $$ 
and
$$ \tilde d = {\rm dim}\, [H^1(\tilde\Gamma\!\sm\!\H^2)] = 1+ [\Gamma(1) : \tilde\Gamma]/6 \,\ge 1+ [\Gamma(1) : \Gamma]/6\, \ge\, d\, . $$ 

  Note that we can have $\,d < \tilde d\,$, as shows the simple example $\,\Gamma= \Gamma(1)$, for which \  $ d = g(\Gamma) = 0\,,\, \nu_\ii(\Gamma) = \nu_\ii(\tilde\Gamma) = g(\tilde\Gamma) =1\, , \, \tilde d = 2\,$. \par \smallskip 

   Furthermore, if $\,F\,$ is a free set of generators of $\,\tilde\Gamma$, and if $\,\gamma\in\Gamma$ represents a generator of $\,\Gamma / \tilde\Gamma$, then $\,\Gamma$ is generated by $\,F\sqcup \{\gamma\}$, with a relation $\,\gamma^n= \phi(F)$ (where $\,n\in\{ 1,2,3,6\}$ is the order of $\,\Gamma / \tilde\Gamma$). Hence the Abelianized $\,{\rm Ab}(\Gamma) := \Gamma / [\Gamma, \Gamma] $ admits the same representation as a free Abelian group, meaning that the Abelianized $\,{\rm Ab}(\tilde\Gamma) := \tilde\Gamma\Big/ [\tilde\Gamma,\tilde\Gamma]$ is isomorphic to a subgroup of $\,\Gamma / [\Gamma, \Gamma] $, such that the quotient be cyclic (of order dividing  $\,n$).  \  
 \if{  
         since \  $ [\Gamma, \Gamma] \subset \tilde\Gamma$ 
$$ {\rm Ab}(\Gamma) := \Gamma / [\Gamma, \Gamma] \equiv \Big(\Gamma / [\tilde\Gamma, \tilde\Gamma] \Big) \Big/ \Big( [\Gamma, \Gamma] / [\tilde\Gamma,\tilde\Gamma] \Big)  \quad \hbox{and} \quad   \Big(\Gamma / [\tilde\Gamma,\tilde\Gamma]\Big) \Big/ \Big( \tilde\Gamma / [\tilde\Gamma, \tilde\Gamma] \Big) \equiv  \Gamma / \tilde\Gamma  $$ 
 $[\Gamma, \Gamma] / [\tilde\Gamma,\tilde\Gamma] = \Big[ \Gamma / [\Gamma, \Gamma] , \Gamma / [\Gamma, \Gamma]\Big]$  must be a subgroup of $\, {\rm Ab}(\tilde\Gamma) := \tilde\Gamma\Big/ [\tilde\Gamma,\tilde\Gamma]$ ; 
}\fi 
Hence, considering the group $\,\widehat{\Gamma}\,$ of additive (real) characters of $\,\Gamma$, we find that : 
$$ \widehat{\Gamma} \equiv \widehat{{\rm Ab}(\Gamma)} \equiv \widehat{{\rm Ab}(\tilde\Gamma)} \equiv \Z^{d'}\, , $$ 
where \ $d'$ is the dimension of the free Abelian group $\,{\rm Ab}(\tilde\Gamma) = \tilde\Gamma\Big/ [\tilde\Gamma,\tilde\Gamma] \equiv \Z^{d'}$, or equivalently, the number of generators of $\,\tilde\Gamma$. \par \medskip 

       Now we have the following important fact, which generalizes Proposition 2.2 of [G-LJ], and, allowing the identification of the cohomology and of the harmonic space of $\,\Gamma\!\sm\!\H^2$,   allows mainly to focus on harmonic forms the asymptotic study of integrals of closed forms. 
\bpro \label{pro.closedform} \quad  $(i)$ \ The number of generators of every free modular group $\,\tilde\Gamma$ equals the dimension of the corresponding harmonic space : \  $ d'=\tilde d := {\rm dim}\, [H^1(\tilde\Gamma\!\sm\!\H^2)]$.  \par
$(ii)$ \  For any modular group $\,\Gamma$, every $C^1$ closed form on $\,\Gamma\!\sm\!\H^2$ is cohomologous to a (unique) harmonic form, element of $\,H^1(\Gamma\!\sm\!\H^2)$. 
\epro
\ub{Proof} \quad  1) \   As ([G-LJ], Prop 2.2), fix a base-point $\,\tilde b\in\tilde\Gamma\!\sm\!\H^2$, $\,b\in \tilde b\,$,  and to any loop $\,\ell\,$ from $\,\tilde b\,$ to $\,\tilde b\,$, associate its lift $\,\bar\ell\,$ in $\,\H^2$ started from $\,b\,$, and $\,\bar\ell(b) = \tilde\ell\, b\in \tilde b\,$. This defines  $\,\tilde\ell\in \tilde\Gamma$  satisfying $\,\widetilde{ ll'} = \tilde l\,\tilde {l'}\,$. \  Then for any $C^1$ closed form  $\,\tilde\omega\,$ on  $\,\tilde\Gamma\!\sm\!\H^2$, set \  ${\ds \tilde\omega(\ell) := \int_\ell \tilde\omega =  \int_{\tilde\pi\circ\bar\ell} \tilde\omega =  \int_{\bar\ell} \tilde\pi_ *\tilde\omega }\,$, \   where  $\,\tilde\pi\,$ denotes the covering projection from $\,\H^2$ onto $\,\tilde\Gamma\!\sm\!\H^2$, which induces a pullback $\,\tilde\pi_ *\,$, mapping $\,\tilde\omega\,$ to a closed form on $\,\H^2$. As $\,\H^2$ is simply connected, $ \tilde\omega(\ell)$ is a function of $\,\bar\ell(b)$, so that we can set \  $ \tilde\omega(\ell) =:  \widehat\omega(\tilde\ell)$, thereby defining an additive character $\,\widehat\omega\,$ on $\,\tilde{\Gamma}$.  \par 
   Now if $\,\widehat\omega=0\,$, meaning that $\,\tilde\omega\,$ is exact, then the primitive $\,{\ds \int_{\tilde b}^{\cdot} \tilde\omega}\, $ is defined on the whole $\,\tilde\Gamma\!\sm\!\H^2$ (since it is arcwise connected), and then constant (as any holomorphic modular function on $\,\Gamma\!\sm\!\H^2$ ; see for example ([L], VI.2.E)), proving that $\,\tilde\omega=0\,$. \  
Hence, the linear map $\,\tilde\omega\mapsto\widehat\omega\,$ is one-to-one from the space of $\,C^1$ closed forms on $\,\tilde\Gamma\!\sm\!\H^2$ into $\,\widehat{\Gamma}$, and, a fortiori,  from $\,H^1(\tilde\Gamma\!\sm\!\H^2)$ into $\,\widehat{\Gamma}$, proving that $\,\tilde d\le d'$. \par \smallskip
   2) \  Reciprocally, notice that the map $\,\ell\mapsto \tilde\ell\,$ defined above induces a morphism $\,\f\,$ from \  $\pi_1(\tilde\Gamma\!\sm\!\H^2)$ into $\,\tilde\Gamma$, since any homotopy $(\ell^s)$ from $\,\ell\,$ to $\,\ell'$ defines a path $(\bar\ell^s(b))$ from $\,\bar\ell(b)$ to $\,\bar\ell'(b)$ in the discrete $\,\tilde\Gamma\,b\,$, forcing $\,\tilde\ell= \tilde\ell'$. \   And if $\,\tilde\ell= 1\,$, then $\,\bar\ell(b)= b\,$ in the simply connected $\,\H^2$, so that $\,\bar\ell\,$ is homotope to the constant  loop $\,b\,$, $\,\ell\,$ is homotope to the constant  loop $\,\tilde b\,$, and thus we have a one-to-one morphism $\,\f\,$ from $\,\pi_1(\tilde\Gamma\!\sm\!\H^2)$ into $\,\tilde\Gamma$. \par 
     Moreover, since $\,\tilde\Gamma$ is discrete and free, it contains only parabolic and loxodromic isometries, and then by Lemma \ref{lem.exp}, any $\,\gamma\in\tilde\Gamma$ can be written  $\,\gamma = \pm\exp(\sigma)$, for a unique $\,\sigma\in sl_2(\R)$. Setting \ $\ell(\gamma)(s) := \pm\exp(s\,\sigma)(\tilde b)$ for $\,0\le s\le 1\,$ and taking the homotopy class $\,\hat\ell(\gamma)\in \pi_1(\tilde\Gamma\!\sm\!\H^2)$ of $\,\ell(\gamma)$, we get a pre-image for $\,\gamma$, with respect to the above morphism $\,\f\,$, which is thus onto. \  By duality ([D], 24.33.2), this yields a one-to-one morphism from $\,\widehat\Gamma$ into $\,H^1(\tilde\Gamma\!\sm\!\H^2)$, hence $\,d'\le \tilde d\,$. \par 
    Hence  $\,\tilde d = d'\,$, and we have exhibited an isomorphism $\,\tilde\omega\mapsto\widehat\omega\,$ from $\,H^1(\tilde\Gamma\!\sm\!\H^2)$ onto $\,\widehat{\Gamma}$. \par \smallskip 
    
    3) \  Fix a basis $(\tilde\omega_1,.., \tilde\omega_{d'})$ of $\,H^1(\tilde\Gamma\!\sm\!\H^2)$, and consider any $C^1$ closed form $\,\tilde\omega\,$ on $\,\tilde\Gamma\!\sm\!\H^2$, which, as described above, defines $\,\widehat\omega\in \widehat{\Gamma}$. By the above, we have reals $\,\alpha_1,..,\alpha_{d'},$ such that \ ${\ds \,\widehat\omega = \sum_{j=1}^{d'} \alpha_j\,  \widehat\omega_j\,}$, implying, as seen above, that \  ${\ds\tilde\omega - \sum_{j=1}^{d'} \alpha_j\,  \tilde\omega_j\,}$ be exact. \par
    So far, we have proved $(i)$ of the statement, and $(ii)$ for the case of a free modular group $\,\tilde\Gamma$. \par \smallskip 
    
    4) \  Let $\,\tilde p\,$ denote the covering projection from $\,\tilde\Gamma\!\sm\!\H^2$ onto $\,\Gamma\!\sm\!\H^2$, which induces a pullback $\,\tilde p_ *\,$, mapping closed smooth differential forms on $\,\Gamma\!\sm\!\H^2$ to closed smooth differential forms on $\,\tilde\Gamma\!\sm\!\H^2$, and harmonic forms on $\,\Gamma\!\sm\!\H^2$ to harmonic  forms on $\,\tilde\Gamma\!\sm\!\H^2$, that is, \ $H^1(\Gamma\!\sm\!\H^2)$ into $H^1(\tilde\Gamma\!\sm\!\H^2)$. \    Note that $\,\tilde p_ *\,$ is necessarily one-to-one (proving again that $\,d\le \tilde d$) : indeed, if a closed smooth form $\,\omega\,$ belongs to its kernel, then we have \ ${\ds  0= \int_{\tilde\ell} \tilde p_ *\omega = \int_{\tilde p\circ \tilde\ell} \omega}$ \   for any loop $\,\tilde\ell\,$  on $\,\tilde\Gamma\!\sm\!\H^2$. \ Now, since the order of the covering group $\,\Gamma/\tilde\Gamma\,$ divides 6,  for any $\,\ell\in H_1(\Gamma\!\sm\!\H^2)$ the lift of $\,\ell^6$ to $\,\tilde\Gamma\!\sm\!\H^2$ is also a loop, meaning that  $\,\ell^6\in \tilde p\circ H_1(\tilde\Gamma\!\sm\!\H^2)$. Hence  \    
${\ds 0 = \int_{\ell^6} \omega = 6 \int_{\ell} \omega}$ \   for any loop $\,\ell\,$  on $\,\Gamma\!\sm\!\H^2$, implying that $\,\omega\,$ must be exact. Considering a primitive of $\,\omega\,$ and using that any holomorphic modular function on $\,\Gamma\!\sm\!\H^2$ is constant (see for example ([L], VI.2.E)), we get $\,\omega =0\,$. \par 
\smallskip 

    5) \  The injectivity of $\,\tilde p_ *\,$ allows to choose the basis $(\tilde\omega_1,.., \tilde\omega_{d'})$ of $\,H^1(\tilde\Gamma\!\sm\!\H^2)$ such that $(\tilde\omega_1,.., \tilde\omega_{d})$ be the image under $\,\tilde p_ *\,$ of a basis $(\omega_1,.., \omega_{d})$ of $\,H^1(\Gamma\!\sm\!\H^2)$ : \  $\tilde\omega_j= \tilde p_ *\omega_j\,$ for $\,1\le j\le d\,$. \ 
    Let us fix a dual basis $(\tilde\ell_1,.., \tilde\ell_{d'})$, that is, a basis of $\,H_1(\tilde\Gamma\!\sm\!\H^2)$ such that \  ${\ds \int_{\tilde\ell_k} \tilde\omega_j = 1_{\{ j=k\}}\,}$ for $\,1\le j,k\le d'$. \par 
    In particular, for $\,1\le j\le d\,$ we have : \  ${\ds 1_{\{ j=k\}} = \int_{\tilde\ell_k} \tilde\omega_j =   \int_{\tilde p\circ \tilde\ell_k} \omega_j\, }$, implying on one hand, for $\,k>d\,$, that $\,\tilde p\circ \tilde\ell_k \equiv 0\,$, and then that  $\,\tilde\ell_k\,$ is a lift to $\,H_1(\tilde\Gamma\!\sm\!\H^2)$ of a loop homotope to 0 in $\,\Gamma\!\sm\!\H^2$, and on the other hand, that $(\tilde\ell_1,.., \tilde\ell_{d})$ is the lift to $\,H_1(\tilde\Gamma\!\sm\!\H^2)$ of a basis $(\ell_1,.., \ell_{d})$ of $\,H_1(\Gamma\!\sm\!\H^2)$. \par \smallskip 
    
    6) \  Consider any $C^1$ closed form $\,\omega\,$ on $\,\Gamma\!\sm\!\H^2$. \   By  3) above, we can write 
$$ \tilde p_*\omega = \sum_{j=1}^{d'} \alpha_j\,\tilde\omega_j + dF \, ,  \quad \hbox{ for some }\; F\in C^2(\tilde \Gamma\!\sm\!\H^2) . $$ 
Integrating this relation along the loop $\,\tilde\ell_k\,$, using 5) above, gives at once $\,\alpha_k=0\,$ for $\,k>d\,$, and \ ${\ds \alpha_k= \int_{\ell_k}\omega\,}$ for $\,1\le k\le d\,$. \par
   Hence \   ${\ds \bar\omega := \omega - \sum_{j=1}^{d} \alpha_j\,\omega_j \,}$ is a $C^1$ closed form $\,\omega\,$ on $\,\Gamma\!\sm\!\H^2$, vanishing on $\,H_1(\Gamma\!\sm\!\H^2)$, and such that \  ${\ds \tilde p_*\bar\omega = dF\,}$. \  For any $\,\tilde b\in\tilde\Gamma\!\sm\!\H^2$, $\,\gamma\in\Gamma$, and any arc $\,c\,$ joining $\,\tilde b\,$ to $\,\gamma\,\tilde b\,$ in $\,\tilde\Gamma\!\sm\!\H^2$, we have : \  ${\ds F(\gamma\,\tilde b) - F(\tilde b) = \int_{c} dF = 
\int_{c} \tilde p_*\bar\omega =  \int_{\tilde p\circ c} \bar\omega = 0\,}$, since $\,\tilde p\circ c\,$ is a loop in $\,\Gamma\!\sm\!\H^2$. \  This proves that $\,F\,$ is $\,\Gamma$-invariant, hence that $\,F=\tilde F\circ \tilde p\,$ for some $\,\tilde F\in C^2(\Gamma\!\sm\!\H^2)$. Finally we have got : \  
$  \tilde p_*(\bar\omega - d\tilde F) = 0\,$, whence \  ${\ds \omega = \sum_{j=1}^{d} \alpha_j\,\omega_j + d\tilde F\,}$,  by 4) above.    $\;\diamond$ 
\par \medskip 

       Moreover, the preceding proposition \ref{pro.closedform} lifts to the modular homogeneous space $\,\Gamma\!\sm\!G\,$ : we can also identifiy the cohomology and the harmonic space of $\,\Gamma\!\sm\!G\,$, and then focus on harmonic forms the asymptotic study of integrals of closed forms on $\,\Gamma\!\sm\!G\,$. 
\bpro \label{pro.closedformG} \   For any modular group $\,\Gamma$, every $C^1$ closed form on $\,\Gamma\!\sm\!G\,$ is cohomologous to a harmonic form, element of $\,H^1(\Gamma\!\sm\!G)$. 
\epro
\ub{Proof} \quad  Recall from Section \ref{sec.fromto} the canonical projection $\,\pi : \Gamma\!\sm\! G \ra \Gamma\!\sm\! \H^2$, whose pullback $\,\pi_*\,$ maps in a one-to-one way the closed forms on $\,\Gamma\!\sm\! \H^2$ to closed forms on $\,\Gamma\!\sm\! G\,$, and the harmonic space $\,H^1(\Gamma\!\sm\! \H^2)$ to the isomorphic subspace $\,\pi_*[H^1(\Gamma\!\sm\! \H^2)]$ of the harmonic space $\,H^1(\Gamma\!\sm\! G)$. \par 
   Fix $\,b= \Gamma\cdot(y,x,\theta)\in \Gamma\!\sm\!G$, and denote by $\,\ell_0\,$ a loop above $\,\pi(b)$, generating $\,H_1(\pi\1(b))$. Fix also a basis $(\ell_1,..,\ell_d)$ of $\,H_1(\Gamma\!\sm\! \H^2)$, dual to the basis $(\omega_1,..,\omega_d)$ of $\,H^1(\Gamma\!\sm\! \H^2)$. Identify   $(\ell_1,..,\ell_d)$ with its lift to $\,H_1(\Gamma\!\sm\! G)$, prescribing merely to each $\,\ell_j\,$  the constant third Iwasawa coordinate $\,\theta\,$.    By Theorem \ref{the.cohom} and by duality (see for example ([D], 24.33.2)), $(\ell_0,\ell_1,..,\ell_d)$ is a basis of $\,H_1(\Gamma\!\sm\!G)$ : indeed, we have \  ${\ds \int_{\ell_0} \pi_*\omega = \int_{\pi\circ\ell_0} \omega = 0\,}$ for any closed form $\,\omega\,$ on $\,\Gamma\!\sm\! \H^2$,  whereas by (\ref{f.omega0}) we have  \  ${\ds \int_{\ell_0} \omega_0 = \int_{\ell_0} d\theta = 2\pi}\,$.  \par
   Consider now any $C^1$ closed form $\,\omega\,$ on $\,\Gamma\!\sm\! G\,$, and set : 
$$  \omega' := \omega - \Big({\ts{1\over 2\pi}} \int_{\ell_0} \omega\Big)\, \omega_0 - \sum_{j=1}^d \Big(\int_{\ell_j} \omega\Big)\, \pi_*\omega_j \, . $$ 
$\,\omega'$ is a $C^1$ closed form on $\,\Gamma\!\sm\! G\,$ such that  $\,{\ds \int_{\ell} \omega = 0\,}$ for any $\,\ell\in H_1(\Gamma\!\sm\!G)$, and then for any loop $\,\ell\,$ on $\,\Gamma\!\sm\!G$. Hence it is exact.  
\if{   
   In the duality between $\,H_1(\Gamma\!\sm\!G)$ and $\,H^1(\Gamma\!\sm\!G)$, $\,\ell_0\,$ is the orthogonal of $\,\pi_*[H^1(\Gamma\!\sm\! \H^2)]$. \par  
   The kernel  of the linear form ${\ds \,\omega\mapsto \int_{\ell_0} \omega}\,$, defined on the vector space of $C^1$ closed forms on $\,\Gamma\!\sm\! G\,$, contains $\,{\cal I}m(\pi_*)$, the image under $\,\pi_*\,$ of the vector space of $C^1$ closed forms on $\,\Gamma\!\sm\! \H^2$. \par 
   Using the metric $((g_{i j}^a))$ of Section \ref{sec.Iwas} and its inverse matrix $((\tilde g_{i j}^a))$, we have the following Petersson-like inner product (recall Formula(\ref{f.Peters})) on the vector space of  closed forms on $\,\Gamma\!\sm\! G$ : 
\beq \label{f.Peters'} \langle \omega , \omega' \rangle  := V(\Gamma\!\sm\! G)\1 \int_{\Gamma\sm G} \tilde g_{i j}^a\,\omega_i \, \omega'_j\,d\mu \, .  \eeq 
$\,\omega \in {\cal I}m(\pi_*)^\bot$ ? 
}\fi 
$\;\diamond$  

\section{Left Brownian motion on ${\ts G= PSL_2(\BR )}$} \label{sec.Brown}  \indf 
   This short section is taken from [F3].  \  Brownian motion $\, g_s = g(z_s,\t_s) = g(y_s, x_s,\t_s)$ on $\,G\,$ has infinitesimal generator $\,\5\,\D^a\,$ and is the left Brownian motion solving the Stratonovitch stochastic differential equation 
$$ d\, g_s = g_s\circ ( \la\, dY_s +\a\, dX_s +\kappa\, a\, dW_s)\, , $$
where $\, (Y_s,X_s,W_s)\,$ denotes a 3-dimensional standard Brownian motion. \par 

   Since a direct calculation shows that 
$$ g( z,\t )\1 dg( z,\t ) = (\sin\t\, dy + \cos\t\, dx)\, {\la\over y} + (\cos\t\, dy - \sin\t\, dx)\,
{\a\over y} +  (y\, d\t + dx)\, {\kappa\over y} \; , $$
we get the differential system 
$$ dy_s = y_s\sin\t_s\circ dY_s + y_s\cos\t_s\circ dX_s = y_s\sin\t_s\, dY_s + y_s\cos\t_s\, dX_s\; , $$  
$$ dx_s = y_s\cos\t_s\circ dY_s - y_s\sin\t_s\circ dX_s = y_s\cos\t_s\, dY_s - y_s\sin\t_s\, dX_s\; , $$
$$ d\t_s = a\, dW_s - \cos\t_s\circ dY_s + \sin\t_s\circ dX_s = a\, dW_s - \cos\t_s dY_s + \sin\t_s\,
dX_s \, .  $$

   Setting \quad 
$ dU_s:= \sin\t_s\, dY_s + \cos\t_s\, dX_s \;$ and $\; dV_s:= \cos\t_s\, dY_s - \sin\t_s\, dX_s \;$, \ 
we get a standard 3-dimensional Brownian motion $\, (U_s,V_s,W_s)\,$ such that 
$$ dy_s = y_s\, dU_s \;\; ,\;\; dx_s = y_s\, dV_s \;\; , \;\; d\t_s = a\, dW_s - dV_s \, .  $$ 

   Hence we see that the projection of our Brownian motion $\, g_s=(y_s,x_s,\t_s)\;$ on the hyperbolic plane $\H^2$, that is to say on the Iwasawa coordinates $\, (y,x)\,$, is simply the standard hyperbolic Brownian motion of $\H^2$, \ and that the angular component $(\t_s)$ is just a real Brownian motion with variance $(1+a^2)$. \ 

\begin{rem}\label{rem.degen} \quad {\rm The degenerate limit-case $a=0\,$ is quite possible for the left Brownian motion $(g_s)\,$. It corresponds to the Carnot degenerate metric on $G$, and to the horizontal left Brownian motion on $G$, associated with the Levi-Civita connexion on $\H^2$.}
\end{rem} 

\section{Asymptotic Brownian windings in $\,\Gamma\!\sm\! G$} \label{sec.exc} \indf 
   Let us denote by 
\beq \label{f.marts}  M^j_t:= \int_{g[0,t]}\omega_j\quad \hbox{ and }\quad  \tilde M^l_t:= \int_{g[0,t]}\tilde\omega_l \eeq 
the martingales obtained by integrating the harmonic forms $\,\omega_j\,$ and $\,\tilde\omega_l\,$ along the paths of the left Brownian motion $(g_s)$. \  
Note that we may as well consider the Brownian motion $\, (g_s)\,$ as living on $G$ or on $\,\Gamma\!\sm\! G\,$. \par 

   Section \ref{sec.Brown} and Theorem \ref{the.cohom} show that 
$$ M^0_t = a\, W_t + \int_0^t \Big( 4\,\Re\,{\ts{\eta'\over \eta}} (x_s,y_s)\, y_s\, dU_s + (4\,\Im\,{\ts{\eta'\over \eta}} (x_s,y_s)\, y_s-1)\, dV_s\Big)\, , $$ 
$$  M^j_t = \int_0^t \Big(\Re f_j(x_s,y_s)\,\, y_s\, dV_s  - \Im f_j(x_s,y_s)\, y_s\, dU_s 
\Big) \quad \hbox{for }\; 1\le j < \nu_\ii(\Gamma) , $$ 
and 
$$ \tilde M^l_t = \int_0^t \Big(\Re \tilde f_l(x_s,y_s)\,\, y_s\, dV_s  - \Im \tilde f_l(x_s,y_s)\, y_s\, dU_s  \Big) \quad \hbox{for }\; 1\le l \le 2{\bf g}(\Gamma) . $$ 

\begin{lem} \label{lem.ctl} \quad The law of \  $t^{-1/2}\,(\tilde M^l_t)_{1\le l \le 2{\bf g}(\Gamma)}\;$ converges towards the centred Gaussian law with diagonal covariance matrix having (on its diagonal) the variances  $\,\langle\tilde\omega_l , \tilde\omega_l\rangle$. 
\end{lem} 
\ub{Proof} \quad By the above, we have some real Brownian motion $(\tilde B^l_s)$ such that 
$$ \tilde M^l_t = \tilde B^l\Big( \langle \tilde M^l\rangle_t\Big) = \tilde B^l\Big( \int_0^t |\tilde\omega_j(z_s)|^2\, ds \Big) , $$ 
and then by scaling, we have the following identity in law (for each $t>0$) : 
$$ t^{-1/2}\, \tilde M^l_t \equiv \tilde B^l\Big( \langle \tilde M^l\rangle_t/t\Big) = \tilde B^l\Big( t\1\int_0^t |\tilde\omega_l(z_s)|^2\,  ds \Big) , $$ 
which by ergodicity converges almost surely to 
$$ \tilde B^l\Big(V(\Gamma\!\sm\!\H^2)\1 \int_{\Gamma\sm\sH^2} |\tilde\omega_l|^2\, dV \Big) = \tilde B^l\Big(\langle\tilde\omega_l , \tilde\omega_l\rangle \Big)\, .  $$ 
Moreover, by ergodicity and by orthogonality of the different $\,\tilde\omega_l\,$, for $\,1\le l < \ell \le 2{\bf g}(\Gamma)$, we have : 
$$ t^{-1}\, \langle \tilde M^l, \tilde M^{l'}\rangle_t  \lra \langle\tilde\omega_l , \tilde\omega_{l'}\rangle = 0\, . $$ 
Hence Knight's Theorem (see ([R-Y], XIII, ¤2, Corollary (2.4))) implies the asymptotic independence of the martingales $(\tilde M^l_t)$. $\; \diamond $
\medskip 

   The following theorem describes the asymptotic Brownian windings in $\,\Gamma\!\sm\! G\,$. 

\begin{theo} \label{the.CV} \   As $\,t\to\ii$, ${\ds \Big( {M^j_t\Big/ t}\, ,\, {\tilde M^l_t\Big/\sqrt{t}\,}
\Big)_{0\le j<\nu_\ii(\Gamma),\, 1\le l\le  2{\bf g}(\Gamma)}\,}$ converges in law towards ${\ds \bigg( \sum_{\ell=1}^{\nu_\ii(\Gamma)} r^\ell_j\, {\mathcal Q}_\ell\, , {\mathcal N}^l \bigg)_{0\le j<\nu_\ii(\Gamma),\, 1\le k\le  2{\bf g}(\Gamma)}}$, where all variables $\,{\mathcal Q}_\ell ,{\mathcal N}^l\,$ are independent, each $\,{\mathcal Q}_\ell\,$ is Cauchy with parameter $\,{h_\ell\over 2\,V(\Gamma\sm\sH^2)}\,$, and $\,{\mathcal N}^l$ is centred Gaussian with variance $\langle\tilde\omega_l , \tilde\omega_l\rangle $. Here $\,h_\ell >0\,$ denotes the width of the solid cusp $\,\CC_\ell$ (already defined in Theorem \ref{the.cohom}). 
\end{theo} 

   Observe the irrelevance of the parameter $a$ in this theorem, which is valid as well in the degenerate case $a=0\,$. \ The reason is that $\, a\,$ was initially the inverse norm of $\LL_\kappa = {\p\over\p\t}\,$, which comes in only in the differential form $\,\omega_0\,$, and which contributes there only to a second order term.
 
\begin{rem}\label{rem.CV}\quad {\rm Theorem \ref{the.CV} is true as well for all finite dimensional marginals : }\parn  
as $t\to\ii$, 
$\;{\ds \Big( {M^j_{c_nt}\Big/ t}\, ,\, {\tilde M^l_{c_nt}\Big/\sqrt{t}\,}
\Big)_{0\le j<\nu_\ii(\Gamma),1\le l\le  2{\bf g}(\Gamma),\, 1\le n\le N}\;}$, for any given $N\in\N^*$ and $\, 0< c_1<..<c_N\,$,  converges jointly towards $\,{\ds \bigg( \sum_{\ell=1}^{\nu_\ii(\Gamma)} r^\ell_j \,{\mathcal Q}^\ell_{c_n}\, , {\mathcal N}^l_{c_n} \bigg)_{0\le j<\nu_\ii(\Gamma),\, 1\le k\le  2{\bf g}(\Gamma), 1\le n\le N}\,}$, where all processes  $\,{\mathcal Q}^\ell , {\mathcal N}^l\,$ are independent, each $\,{\mathcal Q}^\ell\,$ is Cauchy with parameter $\,{h_\ell\over 2\,V(\Gamma\sm\sH^2)}\,$, and $\,{\mathcal N}^l$ is real Brownian with variance $\langle\tilde\omega_l , \tilde\omega_l\rangle $, started from $0$. \end{rem} 

   Note that such statement gives at once the asymptotic law of any finite family of stochastic integrals (along Brownian paths) of harmonic forms, merely by decomposing them in the basis $(\omega_j,\tilde\omega_k)$. \par 
   
  By Proposition \ref{pro.closedformG}, this gives also the asymptotic behavior of  any finite family of stochastic integrals (along Brownian paths) of smooth closed forms, since the normalized contribution of any exact form is clearly negligible in probability (at least for the stationary Brownian motion). 

\section{Geodesics of $\,{\ts G= PSL_2(\BR )}$ and ergodic measures} \label{sec.Geodesics} \indf 
   This section is similar to ([F3], Sections 11,12). 

\subsection{Description of these geodesics} \label{sec.descript}  \indf  
   Recall from Section \ref{sec.Iwas} and Lemma \ref{lem.nat} the metric $((g_{i j}^a))$ (indexed by $\,a\in\R^*$) of $\,G\,$,  as expressed by the Lagrangian $\,L=L (y,x,\t ,\dot y,\dot x,\dot \t)\,$, given in Iwasawa coordinates $\, (y,x,\t )\,$ by :  

\beq \label{f.Lagrang} 2\,L \, =\, y^{-2}\, \dot y^2 + (1+a^{-2})\, y^{-2}\, \dot x^2 + 2a^{-2} y\1\, \dot x \,\dot \t + a^{-2}\, \dot\t^2\, .  \eeq 
The equation of geodesics \  ${ {\p\over\p s} \Big( {\p L\over\p \dot z^j}\Big) =  {\p L\over\p z^j}\;}$ reads here : 
\beq \label{f.eqgeod} (y^{-2}\, \dot y)^{\bf \cdot} = - y^{-3}\, \dot y^2 - (1+a^{-2})\, y^{-3}\, \dot x^2 - a^{-2} y\2\, \dot x \,\dot \t \, ,  \eeq  
and 
\beq \label{f.eqgeod1} (1+a^{-2})\, y^{-2}\, \dot x + a^{-2}\, y^{-1}\, \dot \t = c' \quad\hbox{ and }\quad  y^{-1}\, \dot x + \dot \t = c  \, ,  \eeq  
for two constants $\,c',c\,$. \  Eliminating $\,\dot\t\,$, this gives : 
\beq \label{f.eqgeod2} \dot x = c'y^2 - c\, a\2\,y \quad\hbox{ and }\quad  
(y^{-2}\, \dot y)^{\bf \cdot} = - y^{-3} (\dot y^2 + \dot x^2) -  c\, a\2\,y\2\, \dot x \, ,  \eeq  
Eliminating then $\,c\,a\2\,$ gives : \quad $ (\dot y/y)^{\bf \cdot} = - c'\,\dot x\,$,  \  
or equivalently, for some constant  $\,c''$: 
\beq \label{f.eqgeod5} \dot y/y = c'' - c'\, x \, . \eeq 
Now this implies \  
${\ds  \5\,d \Big({\dot x^2 + \dot y^2 \over y^2}\Big) =  (\dot x/y) c' d y - (\dot y/y) c' d x = 0\,} $,   whence for some non-negative constant $\,C $ : 
\beq \label{f.eqgeod3} (c'\, x-c'')^2 + (c'y - c\, a\2)^2 =\, {\dot x^2 + \dot y^2 \over y^2}\, = C^2\, ,  \quad\hbox{ and }\quad \dot \t = c \, (1+a^{-2}) - c'\, y\, . \eeq  

   In the particular case $\,c'=0\,$, we find \  $c'' x + c\, a\2\, y = c_0$ constant and $\,\dot\theta = c (1+a^{-2})$. \par  
   Hence we see that any geodesic projects on a Euclidian circle or line of $\,\H^2$, and that its projection has constant speed $\,C$ (or energy $\,C^2$). Precisely if $\,c'\not= 0$ :  
\beq \label{f.eqgeod4} c' x = c'' + C \sin\f \,,\, \,\quad c' y = c\, a\2 + C \cos\f \,, \;\hbox{ and } \quad \dot \f = c' y\; , \eeq 
where the last formula results at once from the two preceding ones and from Formula (\ref{f.eqgeod5}). \par 
 If $\,C\not= 0\,$, \   set 
\beq \label{f.eqgeod6} k\, :=\, c\, a\2 C\1  . \eeq 
Then if $\,k\not= 1\,$, according as $\,|k|>1\,$ or $\,|k|<1$, we have :   
$$ \f_s =\,  2\,\arctg\! \bigg[ {\ts\sqrt{k +1\over k-1}}\, \tg\Big( C\,{\ts{\sqrt{k^2-1}\over 2}}\,(s-s_0) \Big)\bigg]  \,\hbox{ or }\;  \f_s =\,  2\,\arctg\! \bigg[ {\ts\sqrt{1+k \over1-k}}\, \th\Big( C\, {\ts{\sqrt{1-k^2} \over 2}}\,(s-s_0)\Big)\bigg] . $$ 

   For $\,|k|>1$, the geodesic projects on a (Euclidian  and hyperbolic) circle totally included in $\,\H^2$, and for $\,|k| <1$, the geodesic projects on a quasi-geodesic of $\,\H^2$ (which is a geodesic if and only if $\,k=0$). In the limiting case $\,|k|=1$, the geodesic projects on an horocycle  of $\,\H^2$, and \  $ \f_s = 2\,\arctg[C(s-s_0)]$. \par 
   Note moreover that by Equations (\ref{f.Lagrang}), (\ref{f.eqgeod1}), (\ref{f.eqgeod3}), and (\ref{f.eqgeod6}), we have :
$$ 2\, L =  {\dot x^2 + \dot y^2 \over y^2} +  a\2 \Big( {\dot x\over y}+ \dot\t\Big)^2 \,= C^2 + a\2 c^2 = (1+a^2 k^2)C^2\, , $$ 
so that prescribing constant speed one to the geodesics implies : 
\beq \label{f.eqgeod7} C\, =\, (1+a^2 k^2)^{-1/2}\,  . \eeq 

   We have in particular established that the energy of any geodesic of $\,\Gamma\!\sm\! G\,$ splits into the constant energy $\,C^2$ of its projection on $\,\Gamma\!\sm\! \H^2$ and the constant energy $\,1-C^2 = a^2 k^2C^2$ of its angular windings about its projection. 

\brem \label{rem.quasigeod} \quad {\rm  The case of main interest for the following is $\,|k| <1$, that is, when the quasi-geodesic $\,\tilde\gamma\,$ of $\,\H^2$ we get intersects $\,\p\H^2$ in two end-points. It is sufficient to consider the case of these two end-points are on the real line, that is, when the quasi-geodesic $\,\tilde\gamma\,$ is a circle, of radius $\,R = C/|c'|$ and centre at height  $\,y_0= kC/c'$. Then the geodesic $\,\psi(\gamma)$ having the same end-points has radius $\,R \sqrt{1-k^2}\,$, and the orthogonal projection $\,p_s\,$ on $\,\psi(\gamma)$ of the point $\,m_s\in \tilde\gamma\,$ having angular coordinate $\,\f_s\,$ has angular coordinate $\,\alpha_s\,$, determined by : \  ${  \sin\alpha_s = { \sqrt{1-k^2}\, \sin\f_s \over 1+ k \cos\f_s}\,}$. As the speed of $\,\tilde\gamma\,$ is $\,C= |\dot m_s| = R |\dot\f_s|/ (y_0+R\cos\f_s) = |\dot\f_s|/ (k+\cos\f_s)$, we find that the geodesic $\,\psi(\gamma)$ must be run at constant speed \  
${ |\dot p_s| = |\dot\alpha_s/\cos\alpha_s| = C \sqrt{1-k^2}\, =  \sqrt{1-k^2\over 1+a^2k^2}\, }$, in order that the distance from $\,m_s\,$ to $\,p_s\,$ remain constant ; necessarily equal to $\,{\rm argch}[({1-k^2})^{-1/2}]$, by standard computation. 
}\erem

  Let us consider the coordinate system $(y,x,\t, u,v,w)$ on $\,T^1G\,$, where $(y,x,\t)$ are the Iwasawa coordinates of the base point $\gamma\in G\,$, and $(u,v,w)$ are the coordinates of the unit tangent vector in the basis $(y{\p\over \p y}, y{\p\over \p x}, {\p\over \p\t})$ of  $\,T^1_{\gamma}G\,$. Thus we have \  $ u^2+v^2+a\2 (v+w)^2 \equiv 1\,$. \par 
  Now, consider the geodesic $(\gamma_s)$ determined by the initial value $(\gamma_0,\gamma_0')\in T^1G\,$ having coordinates $(y_0,x_0,\t_0, u_0,v_0,w_0)$, and let $\,k\,$ be the unique real number determined by the equation \  $v_0+w_0 =  {k\, a^2\over \sqrt{1+a^2k^2}}\,$, and $\,\e := {\rm sign}((1+a^2)v_0+w_0) = {\rm sign}({k\,(1+a^2)\over \sqrt{1+a^2k^2}}-w_0) = {\rm sign}(c')$. \par\smallskip \parn  
  Then by Equations (\ref{f.eqgeod3}), (\ref{f.eqgeod1}), (\ref{f.eqgeod6}), (\ref{f.eqgeod7}),

  the above implies that $(\gamma_s,\gamma_s')$ remains in the leaf $\,L(k,\e)$ having equations : \beq \label{f.feuillek} L(k,\e) := \Big\{ u^2+v^2 = {1\over 1+a^2k^2}\Big\} \bigcap \Big\{ v+w =  {k\, a^2\over \sqrt{1+a^2k^2}}\,\Big\} \bigcap \Big\{ {\rm sign}((1+a^2)v+w) =  \e\Big\} . \eeq

\subsection{Ergodic measures for the geodesic flow on $\,\Gamma\!\sm\! G$} \label{sec.ergmeas} 
\indf  We know by the above section \ref{sec.descript} that any ergodic invariant measure for the geodesic flow on $\,\Gamma\!\sm\! G\,$ must be carried by a leaf $\,L(k,\e)$, for some real $\,k\,$ and $\,\e=\pm 1\,$.  \par \smallskip 

  Note that each geodesic corresponding to the case $\,|k| >1\,$ projects on a periodic curve (circle) in $\,\H^2$, so that there is too few to be said on the asymptotic geodesic behaviour in that case. \  
  Therefore we shall henceforth suppose that  $\,|k| <1$. \par 

\begin{lem} \label{lem.qgeod} \quad  For each $\,k\,$ fixed in $\,]-1,1[\,$ and each $\,\epsilon =\pm1\,$, there is a natural one-to-one map $\,\psi = \psi^{k}_\e\,$ from the leaf $\,L(k,\e)$ (seen as made of geodesics of $\,\Gamma\!\sm\! G\,$ having initial value $\t_0=0$ for their angular part $\t_s$) onto the set of geodesics of $\,\Gamma\!\sm\!\H^2$. \par 
   This map goes as follows : with any geodesic $\,\gamma\,$ of $\,\Gamma\!\sm\! G\,$, associate successively the projection $\,\tilde\gamma\,$ on $\H^2$ of its lift to $G$, and the projection $\psi (\gamma)$ on $\Gamma\!\sm\!\H^2$ of the geodesic of $\H^2$ at bounded distance of $\,\tilde\gamma\,$. \par   
This map makes sense as well at the level of line-elements, and thus defines a homeomorphism from $\,L(k,\e)$ onto $\,T^1(\Gamma\!\sm\!\H^2)\equiv \Gamma\!\sm\! G\,$.  
\end{lem}  
\ub{Proof} \quad  The analysis made in Section \ref{sec.descript} ensures that the map $\,\psi = \psi^{k}_\e\,$ is well defined. Note indeed the necessary $\G$-invariance : if two geodesics $\,\g,\g'$ of $\,G\,$ can be identified modulo some $g\in\G$, then indeed the same $g$ identifies also the geodesics of $\H^2$ at bounded distance of the projections $\,\tilde\gamma , \tilde\gamma'$ of $\,\g,\g'$  on $\H^2$. \par 
   In the reverse direction, to any (oriented) geodesic  $\psi (\g )$ of $\,\Gamma\!\sm\!\H^2$ correspond two quasi-geodesics in $\,\Gamma\!\sm\!\H^2$ at constant distance ${\rm argch}[({1-k^2})^{-1/2}]$, and, owing to the sign $\,\e\,$ which by Formula (\ref{f.feuillek}) and Equation (\ref{f.eqgeod1}) prescribes the sign of $\,c'$ (and then the sign of the height of the centre $\,y_0= kC/c'$ by Remark \ref{rem.quasigeod}), in fact a unique one.  And to this unique quasi-geodesic in $\,\Gamma\!\sm\!\H^2$ is associated by (\ref{f.eqgeod1}) or (\ref{f.eqgeod3}) (for any prescribed initial value $\t_0\,$ of  the angular part) a unique geodesic $\Gamma\gamma$ of $\,\Gamma\!\sm\!G\,$, obviously included in the leaf $\,L(k,\e)$. \    By using furthermore the orthogonal projection in $\H^2$ between our quasi-geodesics and their associated geodesic, we get at once the analogous map at the level of line-elements, with a clear continuity in both directions.   $\;\diamond $  

\begin{rem} \label{rem.feuilles} \quad {\rm Note that in fact each leaf $\,L(k,\e )\,$ splits into 
a continuum of \mbox{sub-leaves :} \  $L(k,\e )=\bigcup_{\t_0\in\sR/2\pi\sZ} L(k,\e,\t_0)$, taking into account the initial value $\t_0$ of the angular part (either at time 0, or above the orthogonal projection of the fixed point $\rt1$ on the quasi-geodesic $\tilde\g$) of the geodesic $\g$. Thus this is indeed the set of its line-elements of each sub-leaf $\,L(C^2,\e,\t_0)$, which is set in one-to-one correspondence with $T^1(\Gamma\!\sm\!\H^2)\equiv \Gamma\!\sm\! G\,$ by the map $\,\psi=\psi_{\e,\t_0}^k$. Note that $\,L(k,\e,\t_0)$ has indeed 3 dimensions, as $G$. \  However, this initial value $\t_0$ will not matter anyway in the following, so that we drop it henceforth, going on with the shorter notation $\,L(k,\e),\, \psi_{\e}^k$. 
}\end{rem} 

\begin{rem} \label{rem.feuilles2} \quad {\rm  According to Section \ref {sec.descript} and Lemma \ref{lem.qgeod} , we have \ 
$$ k^2 = \th^2\Big[\dist\Big(\pi(L(k,\e)),\psi^k_\e(L(k,\e))\Big)\Big]\, , $$
$\,\pi\,$ denoting as in Section \ref{sec.fromto} the canonical projection from $\,\Gamma\!\sm\!G\,$  onto $\,\Gamma\!\sm\!\H^2$. As  \ $\psi^k_\e(L(k,\e))$ is the only geodesic of $\,\Gamma\!\sm\!\H^2$ which is asymptotic to the quasi-geodesic $\,\pi(L(k,\e))$, we see that $\,|k|\,$ is fully determined by the leaf $\,L(k,\e)$, and furthermore, that $\,|k|\,$ is necessarily preserved by any isometry applied to $\,L(k,\e)$. In particular, if $\,0\le k<k'<1\,$ and $\,\e,\e'=\pm 1\,$, then $\,L(k,\e) \cap L(k',\e') = \emptyset\,$. \par 
   As a consequence, note that, counter to the hyperbolic setting,  the geodesic flow on $\,\Gamma\!\sm\! G\,$ is not ergodic (with respect to the normalised Liouville measure of $\,\Gamma\!\sm\! G$) : \  $\bigcup_{0<k<1/2} L(k, \pm1)$ is stable and has measure strictly betwen 0 and 1 (here, as mentioned in Remark \ref{rem.feuilles}, we understand the leaves $\,L(k, \e)$ as 4-dimensional, leaving the initial angular coordinate $\,\t_0\,$ free).  \   
   Observe also that this same set is open in $\,\Gamma\!\sm\!G\,$, proving that the leaves $\,L(k,\e)$ (seen as made of line-elements) are not dense in $\,\Gamma\!\sm\!G\,$. 
}\end{rem} 
\smallskip 

   Now it is known (see [H]) that the Liouville measure on $\,T^1(\Gamma\!\sm\!\H^2)\,$ is invariant and ergodic under the geodesic flow. This fact and Lemma \ref{lem.qgeod} above allow therefore the following. 
\begin{defi} \label{def.Liou} \quad   For any $\,(k, \e )\,$ fixed in $]-1,1[\times\{\pm 1\}$,
denote by $\,\mu^{k}_\e := \mu^\Gamma\circ \psi^{k}_\e\,$, the image of the normalized Liouville measure $\,\mu^\Gamma$ on $\,T^1(\Gamma\!\sm\!\H^2) \equiv \Gamma\!\sm\! G\,$ under the map $\,\psi^{k}_\e\,$ of \mbox{Lemma \ref{lem.qgeod}.} So $\,\mu^{k}_\e\,$ is a probability measure on  the set of line-elements of the leaf $\,L(k,\e )\,$, which is invariant and ergodic under the geodesic flow on $\,\Gamma\!\sm\! G\,$.
\end{defi}
Note that, by Remark \ref{rem.feuilles2}, for $\,|k|<|k'|<1\,$ the measures $\,\mu^{k}_\e , \mu^{k'}_{\e'}\,$ have disjoint supports. 

\section{Asymptotic geodesic windings} \label{sec.Agw} \indf  
   We fix here a leaf $\,L(k,\e )\,$, and endow it with the ergodic invariant probability measure
$\,\mu^{k}_\e\,$ of Definition \ref{def.Liou}. We want to obtain the asymptotic law under $\,\mu^{k}_\e\,$ of 
$$ \Big( t^{-1}\int_{\g [0,t]}\o_j\, ,\, t^{-1/2}\int_{\g [0,t]}\tilde\o_l\,\Big)_{0\le j<\nu_\ii(\Gamma),\, 1\le l\le  2{\bf g}(\Gamma)} \quad \hbox{ as }\;\; t\to\ii\, , $$ 
where the geodesic $\g$ of $\,\Gamma\!\sm\! G\,$ is chosen (at time 0) according to $\,\mu^{k}_\e\,$, and $\,\g [0,t]\,$ denotes this geodesic $\g\,$ run during the time-interval $[0,t]\,$. \par 
   Note that by the $\G$-invariance of the forms $\,\o_j,\tilde\o_l\,$, it makes no difference to think of the geodesics $\,\g\,$ as started in a fundamental domain $D$ and living on $G$, the forms being harmonic on $G$ as well.   \par 
    The following theorem describes the asymptotic geodesic windings in $\,\Gamma\!\sm\! G\,$, under the ergodic measures of Section \ref{sec.ergmeas}. A minor mistake in [F3], concerning the contribution of $\,d\theta\,$,  is corrected here. 

\begin{theo} \label{the.geo}\quad Let us consider a fixed leaf $\,L(k,\e )\,$ (defined in Section 
\ref{sec.descript}) of $\,\Gamma\!\sm\! G\,$, with $\,|k| <1\,$, endowed with the ergodic invariant probability measure $\,\mu^{k}_\e$ of \mbox{Definition \ref{def.Liou}.}  \par  
Then the law under \  $\mu^{k}_\e = \mu^{k}_\e (d\g )$ \  of \quad 
${\ds \Big( t^{-1}\int_{\g [0,t]}\o_j\, ,\, t^{-1/2}\int_{\g [0,t]}\tilde\o_l\,\Big)_{0\le j<\nu_\ii(\Gamma),\, 1\le l\le  2{\bf g}(\Gamma)} } $ \parn   converges as $t\to\ii\,$ to the law of \quad 
$${\ds \left(  {(1+a^2) k \over \sqrt{1+a^2k^2}} \, 1_{\{ j=0\}} + 2 \Big({1-k^2\over 1+a^2k^2}\Big)^{1/2} \sum_{\ell=1}^{\nu_\ii(\Gamma)} r^\ell_j\, {\mathcal Q}_\ell \, , \Big({4(1-k^2) \over 1+a^2k^2}\Big)^{1/4} {\mathcal N}^l \right)_{0\le j<\nu_\ii(\Gamma),\, 1\le l\le  2{\bf g}(\Gamma)}}\, ,  $$  
where all variables $\,{\mathcal Q}_\ell ,{\mathcal N}^l\,$ are independent, each $\,{\mathcal Q}_\ell\,$ is Cauchy with parameter $\,{h_\ell\over 2\,V(\Gamma\sm\sH^2)}\,$, and $\,{\mathcal N}^l$ is centred Gaussian with variance $\langle\tilde\omega_l , \tilde\omega_l\rangle $. Here $\,h_\ell >0\,$ denotes the width of the solid cusp $\,\CC_\ell$ (defined in Theorem \ref{the.cohom}). 
\end{theo} 

   Note a clear difference between the Brownian and geodesic behaviors : mainly, here (counter to the Brownian case) the $\,d\t$-part of the form $\,\o_0\,$ is responsible for a non-negligible asymptotic contribution.   Moreover the parameter $a$ now appears in the limit law. \parn  
This makes a noteworthy contrast with the hyperbolic case (see [E-F-LJ1], [E-F-LJ2], [F2]).  \par  
    This difference appears in Lemma \ref{lem.w12} below, whereas once the $\,d\t$-part has been moved away, the remaining asymptotic law is essentially the same as the Brownian one, given by Theorem \ref{the.CV}.  So the remaining task, following [F3], will be then (in Section \ref{sec.proofg} below) mainly to compare on $\H^2$ the geodesic paths to the Brownian paths, somewhat in the spirit of the methods already employed in [E-F-LJ1], [E-LJ], [F2], [LJ2], but in a simpler way, taking advantage of the harmonicity of the forms $\o_j,\tilde\o_l\,$, somewhat as in [LJ1]. 
\par \medskip   
       The following lemma (in which the minor mistake in [F3] concerning the contribution of $\,d\theta\,$ is corrected) reduces the study along the geodesics of $G$ to a study along the geodesics of $\H^2$. 
\begin{lem}\label{lem.w12} \   The asymptotic law of  
${\ds \Big( t^{-1}\int_{\g [0,t]}\o_j\, ,\, t^{-1/2}\int_{\g [0,t]}\tilde\omega_l\,\Big)_{0\le j<\nu_\ii(\Gamma),\, 1\le l\le  2{\bf g}(\Gamma)} } $ under $\mu^{k}_\e (d\g)$  is the same as the asymptotic law under the Liouville measure $\mu^\Gamma$ on $T^1(\Gamma\!\sm\!\H^2)$ of : 
$$ \left(  {\ts{(1+a^2) k \over \sqrt{1+a^2k^2}}}\,1_{\{j=0\}}  + {\ts\Big({1-k^2\over 1+a^2k^2}\Big)}^{1/2} t^{-1} \int_{\g [0,t]}\o_j\, , {\ts\Big({1-k^2\over 1+a^2k^2}\Big)}^{1/4} t^{-1/2} \int_{\g [0,t]}\tilde\omega_l \right)_{0\le j<\nu_\ii(\Gamma),\, 1\le l\le  2{\bf g}(\Gamma)} . $$ 
\end{lem}  
\ub{Proof} \quad Let us deal first with $\,\tilde\omega_l = \pi^*\tilde\omega_l\,$. By Definition \ref{def.Liou} and Lemma \ref{lem.qgeod}, we just have to compare \  
${\ds t^{-1/2}\int_{\g [0,t]}\tilde\omega_l = t^{-1/2}\int_{\tilde\g[0,t]}\tilde\omega_l }$ \  with \  
${\ds t^{-1/2}\int_{\psi (\g )[0,t]}\tilde\omega_l \,}$. \parn  
Now use that on $\H^2$, \   $\tilde\omega_l = dF_l\,$ is exact, and 
recall from Remark \ref{rem.quasigeod} that 
$$\dist\Big( F_j[\tilde\g(t)] , F_j\Big[\psi (\g )({\ts\sqrt{1-k^2\over 1+a^2k^2}}\,t)\Big]\Big) = (1-k^2)^{-1/2}\, , $$
to get :  
$$ \left|\int_{\tilde\g[0,t]} \tilde\omega_l  - \int_{\psi (\g )\left[ 0, \left( {1-k^2\over 1+a^2k^2}\right)^{1/2}\,t \right]} \tilde\omega_l  \right| = \Big| F_j[\tilde\g(t)] - F_j\Big[\psi (\g )({\ts\sqrt{1-k^2\over 1+a^2k^2}}\,t)\Big] - F_j[\tilde\g(0)] + F_j[\psi (\g )(0)]\, \Big| $$ 
$$ \le\, 2\,\| \tilde\omega_l \|_\ii\Big/\sqrt{1-k^2}\; . $$ 
This shows that \  
${\ds t^{-1/2}\left(\int_{\tilde\g[0,t]} \tilde\omega_l  -\int_{\psi (\g )\left[ 0, \left( {1-k^2\over 1+a^2k^2}\right)^{1/2}\,t \right]} \tilde\omega_l \right)}$ \  
goes uniformly  to 0, proving the result relating to the regular forms $\,\tilde\omega_l\,$. \par 
   Now we have to deal with the singular forms $\,\omega_j\,$, which from Theorem \ref{the.cohom} equal $\,\pi^*\omega'_j\,$ for $\,1\le j<\nu_\ii(\Gamma)$, whereas  $\,\omega_0 = d\t + \pi^*\omega'_0\,$.  Thus we can handle the $\,\omega'_j\,$ as the $\,\tilde\omega_l\,$ above, using Theorem \ref{the.cohom} again, to get : 
$$ \int_{\g [0,t]}\omega'_j = \int_{\tilde\g[0,t]}\omega'_0 = \int_{\psi (\g )\left[ 0, \left( {1-k^2\over 1+a^2k^2}\right)^{1/2}\,t \right]}\omega'_j + \O\Big(\tilde{Y}(\psi\circ\g ,0) + \tilde{Y}(\psi\circ\g ,t)\Big)\, , $$  
where \quad $\tilde{Y}(g,t) := \max_{1\le \ell\le \nu_\ii(\Gamma)}\limits\, \sup\Big\{ \yt_\ell (z)\,\Big|\, \dist\Big( g\Big[{\ts\sqrt{1-k^2\over 1+a^2k^2}}\,t\Big],z\Big)\le 1/\sqrt{1-k^2}\,\Big\}$. \parn 
On the other hand we have by Equations (\ref{f.eqgeod3}), (\ref{f.eqgeod6}), and (\ref{f.eqgeod7}) : 
$$ t\1\int_{\g [0,t]} d\t = c\,(1+a\2) - t\1\int _{\g [0,t]}  c'y_s\,ds\, \lra\, c\,(1+a\2) - c'\lim_t y(\gamma(t)) $$ 
$$ = c\,(1+a\2) = {(1+a^2) k\Big/\sqrt{1+a^2k^2}}\, . $$ 

   Therefore, by Definition \ref{def.Liou}, the asymptotic law of \  ${\ds t\1\int_{\g [0,t]}\omega_j}$ \   under $\,\mu^{k}_\e (d\g )\,$ is the same as the asymptotic law under $\,\mu(dg)$ of  :  
$$ {(1+a^2) k\over\sqrt{1+a^2k^2}}\,1_{\{j=0\}} + t\1\int_{g \left[ 0, \left( {1-k^2\over 1+a^2k^2}\right)^{1/2}\,t \right]}\omega'_j + \O \Big(\tilde{Y}(g ,0) + \tilde{Y}(g ,t)\Big)\Big/ t \, . $$ 
Observe further that under $\,\mu^\Gamma(dg)$ the process $\,\tilde{Y}(g ,t)$ is stationary, so that the last term above asymptotically vanishes in probability as $\,t\to\ii\,$. \ Hence we have shown that \  the asymptotic law of \  ${\ds t\1\int_{\g [0,t]}\omega_j}$ \  under $\,\mu^{k}_\e (d\g)$ is the same
as the asymptotic law under $\,\mu (dg)$ of
$$ {(1+a^2) k\over\sqrt{1+a^2k^2}}\,1_{\{j=0\}} + t\1\int_{g \left[ 0, \left( {1-k^2\over 1+a^2k^2}\right)^{1/2}\,t \right]} \omega'_j  \,  . $$ 
Finally the result is valid jointly for the terms with singular forms $\,\omega_j\,$ and regular forms $\,\tilde\omega_l\,$, since for each the negligible contributions vanish in probability. It remains only to replace $\,t\,$ by $\, \Big({1-k^2\over 1+a^2k^2}\Big)^{-1/2}\, t\,$.  $\;\diamond $ \par 

\section{Equirepartition in $\,\Gamma\!\sm\! G\,$ of large geodesic spheres} \label{sec.Erlgs} \indf  
    Recall that $\,\mu\,$ (defined in Section \ref{sec.Geom} by \  ${ \mu (dg) := { dx\, dy\, d\t\over 4\pi^2\, y^2}}\,$) denotes the Liouville measure on $\,G\,$, and that $\,\mu^\Gamma$ denotes its normalized projection on $\,\Gamma\!\sm\! G\,$, identified with a law on left $\,\Gamma$-invariant functions  on $\,G$ ; and ${\rm covol}(\Gamma)$ denotes the finite hyperbolic volume of $\,\Gamma\!\sm\! \H^2$.   \par \smallskip 
   Fix some $\,g\in G\,$, and a PoincarŽ half-plane model $\,\P^2$ for $\,\H^2$ such that $\,\pi(g) = e_0=(0,1)$. \   For any $\, R>r>0\,$, let \  $B_{r} :=B(e_0 ;r)$ denote the geodesic ball of radius $\,r\,$ in $\,\P^2$, centred at $\,e_0=(0,1)$, and (denoting by $\Theta$ the geodesic flow) consider : 
$$ T_R(r) := (T^1\! B_r)\, \Theta_{R} = PSO(2)\,\Theta_{[0,r]}\,PSO(2)\, \Theta_{R} \,\subset\,  T^1(B_{R+r}\moins B_{R-r})\subset T^1\P^2\equiv G\, . $$
   Set then  :   
$$  \f_R := \sum_{\gamma\in\Gamma} 1_{\gamma\,T_R(r)} = \f_0\circ\Theta_{R}\, , \quad \hbox{so that } \; \int\f_R\,d\mu^\Gamma \,= \,{{\rm vol}(B_r)\over {\rm covol}(\Gamma)}\, = \,{2\pi\,(\ch r-1)\over {\rm covol}(\Gamma)}\, \in\, ]0, \ii [ . $$  
Indeed, by definition of $\,\mu , \mu^\Gamma$, considering some fundamental domain $\,D$ of $\,\Gamma$,  we have : 
$$ \int\f_R\,d\mu^\Gamma = \int\f_0\,d\mu^\Gamma = \sum_{\gamma\in\Gamma} \mu^\Gamma (T^1D\cap \gamma\, T^1B_r) = {\rm covol}(\Gamma)\1\sum_{\gamma\in\Gamma} {\rm vol}(D\cap \gamma\, B_r) $$ 
$$ = {\rm covol}(\Gamma)\1\sum_{\gamma\in\Gamma} {\rm vol}(B_r\cap\gamma\,D) = \,{{\rm vol}(B_r)\over {\rm covol}(\Gamma)}\, = \,{\int_0^{2\pi}\!\int_0^r \sh\rr\,d\rr\,d\t \over {\rm covol}(\Gamma)}\, =  \,{2\pi\,(\ch r-1)\over {\rm covol}(\Gamma)}\, . $$ 
Similarly, for any $\,f\in L^2(\mu^\Gamma)$ we have : 
$$ {\ts{{\rm covol}(\Gamma)\over 2\pi\,|a|}}\times  \int  f\times \f_R\,\, d\mu^\Gamma =  \int_{T^1D} \sum_{\gamma\in\Gamma}  1_{\gamma\,T_R(r)}\, f\,d\mu =  \sum_{\gamma\in\Gamma}  \int 1_{T_R(r)\cap\gamma T^1D}\, f\, d\mu =  \int_{T_R(r)}  f\, d\mu \, . $$ 
Moreover, for any $\,r>0\,$, by the discontinuity of the action of $\,\Gamma\,$ on $\,\H^2$, $\,\f_0\in L^\ii(\mu^\Gamma)$. \par  

   Applying now the mixing theorem, for any $\,f\in L^2(\mu^\Gamma)$ and $\,r>0$, we have as $\,R\to\ii$ : 
$$  \int \f_R\times f\,\, d\mu^\Gamma  =   \int \f_0\circ\Theta_{R}\times f\,\, d\mu^\Gamma \lra  \int \f_0\,d\mu^\Gamma\times \int f\,d\mu^\Gamma \, , $$ 
or equivalently : 
$$ {\ts{2\pi\,|a|\over {\rm vol}(B_r)}} \times \int_{T_R(r)}  f\, d\mu \,=\, {\ts{{\rm covol}(\Gamma)\over {\rm vol}(B_r)}}\times  \int  f\times \f_R\,\, d\mu^\Gamma \lra   \int f\,d\mu^\Gamma \, . $$ 

    Therefore, for any $\,r>0$, the probability law $\,\mu^\Gamma$ on $\,\Gamma\!\sm\! G\,$ is the weak limit, as $\,R\to\ii$, of the projection on $\,\Gamma\!\sm\! G\,$ of the normalized restrictions $\, {2\pi\,|a|\over {\rm vol}(B_r)}\, 1_{T_R(r)} \, \mu\,$ of $\,\mu\,$ to the large shells $\,T_R(r)$.   
\par \smallskip 

   Suppose then that the function $\,f\,$ on $\,\Gamma\!\sm\! G\,$ is continuous and compactly supported, and then uniformly continuous as a (left $\,\Gamma$-invariant) function on $\,G=PSL_2(\R)$. 
   
   Then, it is easily seen (and verified by a standard computation) that the maximal angular deviation in a section of diameter $\,r\,$ of a large thin shell $\,T_R(r)$, between geodesics arriving from $\,T^1\! B_r\,$, is equivalent to $\,{2Rr\over R-1}\,$, as $\,r\sea 0\,$ and $\,R\ge 2\,$. This implies 
that 
$$ \xi\in T^1\! B_r\,,\; \rr\in PSO(2)\, , \,\dist_{\sH^2}[\pi(\xi\Theta_R), \pi(\rr\Theta_R)] <r \,\Ra\, 
\dist_{PSL_2(\sR)}(\xi\Theta_R\, ,\,\rr\Theta_R) = \O(r)\, . $$ 
  
   Hence, denoting by $\,d\rr\,$ the normalized uniform measure on $PSO(2)$, we have :  
$$  \Big|\, {\ts{2\pi\,|a|\over {\rm vol}(B_r)}}\! \int_{T_R(r)}  f\, d\mu - \int_{PSO(2)} f(\rr\, \Theta_{R})\, d\rr \,\Big| \lra  0 \quad \hbox{ as }  \;  r\sea 0\, ,  $$ 
uniformly with respect to $\,R\ge 2\,$. \par 
   This means that  the mean of $\,f\,$ on the large thin shell $\,T_R(r)$ converges, as $\,r\sea 0\,$ and $\,R\ge 2\,$, to its mean on the geodesic sphere \  $S_{R} := (T^1_{e_0}\P^2)\Theta_{R}\equiv PSO(2) \Theta_{R} \subset T^1(S(e_0 ; R))$.  \par 
   Finally, we get that :  
$$  \Big|  \int f\,d\mu^\Gamma  - \int_{S_R} f \,\Big|\, \le\, \Big|  \int f\,d\mu^\Gamma - {\ts{2\pi\,|a|\over {\rm vol}(B_r)}}\!\int_{T_R(r)}  f\, d\mu \,\Big| +\Big|\, {\ts{2\pi\,|a|\over {\rm vol}(B_r)}}\! \int_{T_R(r)}  f\, d\mu - \int_{S_R} f \,\Big|   $$ 
is arbitrary small for $\,r\,$ fixed so small that the first term on the right hand side be small enough, and then for large enough $\,R\,$.  \   Forgetting the irrelevant choice of the PoincarŽ model $\,\P^2$ for $\,\H^2$ and of its base point $\,e_0\,$, and projecting on $\,\Gamma\sm G$, we can see $\,S_R\,$ as the sphere of radius $R\,$ in $\,\Gamma\!\sm\! G\,$, centred at $\,\Gamma g\in \Gamma\!\sm\! G\,$. \par  \smallskip 
   This proves the following. 
\bthe \label{the.Erlgs} \quad  The normalized Liouville measure $\,\mu^\Gamma$ on $\,T^1(\Gamma\!\sm\!\H^2)\equiv \Gamma\!\sm\! G\,$ is the weak limit as $\,R\to\ii\,$ of the uniform law on the geodesic sphere  $\,\Gamma g\,PSO(2) \Theta_{R}$ of $\,\Gamma\!\sm\! G\,$ having radius $\,R\,$ and fixed center $\,\Gamma g\in \Gamma\!\sm\! G$ : for any  compactly supported continuous function $\,f\,$ on $\,\Gamma\!\sm\! G\,$, denoting by $\,d\rr\,$ the uniform law on $PSO(2)$, we have 
$$  \int f\,d\mu^\Gamma = \lim_{R\to\ii} \int_{PSO(2)} f(\Gamma g\,\rr \,\Theta_{R})\, d\rr \, . $$ 
\ethe       
\brem \label{rem.Erlgs} \   {\rm  This result and its proof are contained in [E-MM]. Such equidistribution result goes back to [R].  
}\erem 

   For any compactly supported continuous function $\,h\,$ on $\,L(k,\e)$ (recall Definition \ref{def.Liou}), applying Theorem \ref{the.Erlgs} to $\,h\circ(\psi^{k}_\e)\1\,$ (this is licit according to Lemma \ref{lem.qgeod}), we get : 
$$ \int_{L(k,\e)} h\, d\mu^{k}_\e = \int_{\Gamma\!\sm\! G} h\circ(\psi^{k}_\e)\1\, d\mu^{\Gamma} = 
\lim_{R\to\ii} \int_{PSO(2)} h[(\psi^{k}_\e)\1(\Gamma g\,\rr \,\Theta_{R})]\, d\rr  \, , $$ 
hence the following equidistribution result, reminiscent of a multi-dimensional ergodic theorem. 
\bcor \label{cor.Erlgs} \quad  For any $(k, \e )$ fixed in $]-1,1[\times\{\pm 1\}$ and any $\,g\in G$, the probability measure $\,\mu^{k}_\e\,$ on the leaf $\,L(k,\e)$ (recall Definition \ref{def.Liou}) is the weak limit as $\,R\to\ii\,$ of the uniform law on the geodesic quasi-sphere $(\psi^{k}_\e)\1(\Gamma g\,PSO(2) \Theta_{R})$. 
\ecor


\section{Synthetic proof of Theorem \ref{the.CV}}  \label{sec.SyntPr} \indf 
    The proof of Theorem 2 in ([F3], Section 10) essentially applies here, with minor modifications. Thus this  section presents a somewhat sketched proof of Theorem \ref{the.CV}, containing all ingredients, but not all details, for which we refer to [F3]. \par \smallskip 
    
    The slow windings (about the handles) $\tilde M^l_t$ were already easily handled in Lemma \ref{lem.ctl}. Hence we must now deal with the singular windings (about the cusps) $M^j_t$, and then establish the asymptotic independence  of both types, which is the main difficulty of the whole proof and is widely responsible for its length. \par\smallskip 
    
    1) \   To proceed, we first cut the solid cusps at some high level $\,r>0$, considering (for $0\le j<\nu_\ii(\Gamma)$ and $1\le \ell \le \nu_\ii(\Gamma)$) the martingales 
\beq \label{f.Mjlr}  M^{j,\ell,r}_t := 1_{\{ r^\ell_j\not= 0\}}\, (r^\ell_j)\1 \int_0^t 1_{\{ \yt_\ell(s)>r\}}\, dM^{j}_s \, , \; \hbox{ where } \; \yt_\ell(s) := 1_{\{ g_s\in\CC_\ell\}}\,\yt_\ell(g_s) , \eeq 
$\,\yt_\ell\,$ (defined in Section \ref{sec.geomG}) being the height in the cusp $\,\CC_\ell\,$. \par  

   Observe that the martingale ${\ds \,M^{j}_t-\sum_{\ell =1}^{\nu_\ii(\Gamma)} r^\ell_j\, M^{j,\ell,r}_t}$, locally constant out of the compact \parn ${\ds \bigcap_{\ell =1}^{\nu_\ii(\Gamma)}\{\yt_\ell \le r\}}$, has bounded quadratic variation, so that \   
${\ds \Big( M^{j}_t-\sum_{\ell =1}^{\nu_\ii(\Gamma)} r^\ell_j\, M^{j,\ell,r}_t}\Big)\Big/\sqrt{t}$  \  converges in law and \  
${\ds \Big( M^{j}_t-\sum_{\ell =1}^{\nu_\ii(\Gamma)} r^\ell_j\, M^{j,\ell,r}_t}\Big)\Big/ {t}$ \   goes to 0 in $L^2$-norm, as $\,t\to\ii\,$. \  Set \ 
\beq \label{f.martMlr} M^{\ell,r}_t := \int_0^t 1_{\{ \yt_\ell(s)>r\}}\, d\xt_\ell(s) \,, \; \hbox{ where } \; \,\xt_\ell(s) := 1_{\{ g_s\in\CC_\ell\}}\,\xt_\ell(g_s) .  \eeq    
Owing to Theorem \ref{the.cohom} and Section \ref{sec.Brown}, we see that \   
${\ds \Big( M^{j,\ell,r}_t - M^{\ell,r}_t\Big)\Big/ {t} = t\1\int_0^t\O(1)\, dV^\ell_s\,}$ \   goes also to 0 in $L^2$-norm. \par 
   We have therefore only to study the martingales $(M^{\ell,r}_t)$, instead of the $(M^{j}_t)$. \par \smallskip 
   2) \   Consider then a discretization of the excursions of the Brownian motion $(g_t)$ in the cusps :
it enters the shortened solid cusp $\{\yt_\ell > r+\sqrt{r}\}$ and exits $\{\yt_\ell \ge r\}$ for the $n$th time within the interval of time say $[\tau_n^\ell,\sigma^\ell_n]$, during which it performs an elementary winding \
\beq \label{f.philr} \f_n^\ell = \f_n^\ell(r) := \int_{\tau_n^\ell}^{\sigma^\ell_n} d\xt_\ell(s). \eeq 
Depending only on planar hyperbolic Brownian motions (recall Section \ref{sec.Brown}), these elementary windings are independent, and independent from the points on the level $\{\yt_\ell = r\}$ at which the excursions start, and are easily (and classically) seen to have a Cauchy law, of parameter $\,\sqrt{r}\,$.  A random number $\,\la^\ell_t = \la^\ell_t(r)\,$ of these windings is performed till time $\,t\,$. By ergodicity, we have $\,\lim_{t\to\ii}\limits \la^\ell_t/t =:\rr^\ell_r\,$ almost surely. \  Otherwise the Markov property implies the independence of the excursion durations $\{\sigma^\ell_n-\tau^\ell_n\,|\,n\in\N^*\}$, so that by the law of large numbers \  ${\ds N\1\sum_{n=1}^N (\sigma^\ell_n-\tau^\ell_n)}\,$ goes almost surely to $\,\E(\sigma^\ell_1-\tau^\ell_1)$ as $N\to\ii$. \  By an obvious comparison between ${\ds \,\sum_{n=1}^{\la^\ell_t } (\sigma^\ell_n-\tau^\ell_n)}\,$ and ${\ds \,\int_0^t 1_{\{ \yt_\ell(s)>r\}}\, ds\,, \int_0^t 1_{\{ \yt_\ell(s)>r+\sqrt{r}\}}\, ds\,}$, and by the ergodic theorem, we deduce that 
$$ { V(\{ \yt_\ell >r+\sqrt{r}\}) \over V(\Gamma\sm\H^2)} \le \rr^\ell_r\times \E(\sigma^\ell_1-\tau^\ell_1) \le { V(\{ \yt_\ell >r\}) \over V(\Gamma\sm\H^2)} \, . $$ 
Now, on one hand it is easily computed that $\,\E(\sigma^\ell_1-\tau^\ell_1)= 2\,\log(1+r^{-1/2})$, and on the other hand, we have : \  ${\ds V(\{ \yt_\ell >r\}) = \int_{[0,h_\ell]\times]r,\ii[} \tilde y_\ell\2\,{d\tilde x_\ell \, d\tilde y_\ell} = h_\ell/r\,}$, by definition of the width $\,h_\ell\,$ (recall Theorem \ref{the.cohom}). \  Hence we find that \  
\beq \label{f.nbreexc}  \lim_{r\to\ii}\, \sqrt{r}\, \lim_{t\to\ii}\limits {\la^\ell_t(r)\over t}  \,=\, \lim_{r\to\ii}\, \sqrt{r}\, \rr^\ell_r \,=\, { h_\ell  \over 2\,V(\Gamma\sm\H^2)}\,\, . \eeq  

   3) \  Let us now analyse further the behaviour of the martingales $\,M^{\ell,r}_t\,$ of formula (\ref{f.martMlr}), by means of the above excursions. There are possibly two incomplete excursions, namely the very first one, the winding contribution (divided by the normalisation $\,t$) of which almost surely vanishes, and the very last one, which exists only when the Brownian motion at time $\,t\,$ visits the solid cusp $\{ \yt_\ell >r\}$, which is the case only with probability $\,\O(1/r)$, so that 
its winding contribution (letting $\,r\to\ii$) eventually vanishes in probability. \  Hence the only non-negligible contribution of the martingale $\,M^{\ell,r}_t\,$ comes from \  
${\ds \sum_{n=1}^{\la^\ell_t } \f_n^\ell }\,$. \   Then, using again that $\,\lim_{t\to\ii}\limits \la^\ell_t/t =\rr^\ell_r\,$, taking advantage of the above observation that ${\ds \sum_{n=1}^{N} \f_n^\ell }\,$ constitutes a discretized Cauchy process (of parameter $\,\sqrt{r}\,$), and using the scaling property and the right continuity of a Cauchy process, we see that \  ${\ds t\1\sum_{n=1}^{\la^\ell_t} \f_n^\ell }\,$, hence $\,M^{\ell,r}_t/t\,$, has, in probability, the same asymptotic behaviour as \  ${\ds t\1\sum_{n=1}^{[\rr^\ell_r\,t] } \f_n^\ell }\,$ ; and as $\,t\to\ii$, this last process converges in law towards a Cauchy variable of parameter $\,\sqrt{r}\,\rr^\ell_r\,$. \par \smallskip 

   4) \   To establish the asymptotic independence of Theorem \ref{the.CV}, we need to approach also the slow windings martingales  $\,\tilde M^k_t\,$ (recall Formula \ref{f.marts}), by martingales that are supported in the complement $\,K_r\,$ of all solid cusps $\{ \yt_\ell >r+\sqrt{r}\,\}$, in order to be able to take advantage of the Markov property, from which independence can then derive. \par 
   Precisely, let us order all stopping times $\{\tau_n^\ell\,,\sigma^\ell_n\,|\,1\le \ell\le \nu_\ii(\Gamma), \, n\in\N^*\}$ into a unique strictly increasing sequence $\,( ..< \tau_n<\sigma_n< ..)\,$, fix any $\,\tilde q\in\R^{2{\bf g}(\Gamma)}$, and consider \  ${\ds  J_n^q := \sum_{k=1}^{2{\bf g}(\Gamma)} \tilde q_k \int_{\sigma_n}^{\tau_{n+1}} d\tilde M^k_t \;}$.  \   Let \  ${\ds \la_t := \sum_{\ell =1}^{\nu_\ii(\Gamma)} \la^\ell_t \,}$ be the total number of excursions performed till time $\,t\,$. \   The same argument as for Lemma \ref{lem.ctl} proves that,  as $\,t\to\ii\,$,  \parn  
${\ds \sum_{k=1}^{2{\bf g}(\Gamma)} \tilde q_k\,t^{-1/2}\, \tilde M^k_t -  t^{-1/2}\sum_{n=1}^{\la_t} J_n^q  \;}$ is asymptotically  $\,\O(1/r)$, provided we can handle the last excursion in the compact core $\,K_r\,$, alive at time $\,t$ ; now, considering the quadratic variation and using the integrability of $(\tau_{n+1}-\sigma_n)^2$ and that $\,\la_t/t\,$ is bounded in probability, it is easily seen that this last excursion in $\,K_r\,$ has a contribution which  vanishes in probability. \par 
   Hence we can asymptotically substitute \  ${\ds \sum_{n=1}^{\la_t} J_n^q\;}$ for the martingale \  ${\ds \sum_{k=1}^{2{\bf g}(\Gamma)} \tilde q_k\,\tilde M^k_t\,}$. \par 
   Consider now the Markov chain $(Z_{\sigma_n},Z_{\tau_{n+1}})$ induced, for any fixed $\,r\,$, by the Brownian motion $(Z_t)$ on $\,H^2/\Gamma\,$, which is known to be stationary and ergodic under the so-called Palm probability measure $\,\chi\,$ induced by the volume measure on the union of all boundaries $\{ \yt_\ell =r\}$, $\{ \yt_\ell =r+\sqrt{r}\,\}$ of the solid cusps. The transition operator of this induced Markov chain has a sprectral gap in $\,L^2(\chi)$, which implies that correlations between durations $(\tau_{n+1}-\sigma_n)$ decay exponentially fast. \par  
   This implies in turn that (the quadratic variation of) \  ${\ds \sum_{n=1}^{\la_t} J_n^q - \sum_{n=1}^{[\rr_r\,t]} J_n^q\;}$ goes to 0 in probability, where \  ${\ds \rr_r := \sum_{\ell =1}^{\nu_\ii(\Gamma)} \rr^\ell_r\,}$ is deterministic : we can substitute \  ${\ds \sum_{n=1}^{[\rr_r\,t]} J_n^q \;}$ for  \  ${\ds \sum_{n=1}^{\la_t} J_n^q\;}$.  \par  \smallskip 
   
   5) \   Consider any $(q,\tilde q)\in \R^{\nu_\ii(\Gamma)}\times\R^{2{\bf g}(\Gamma)}$, and 
\beq \label{f.loilim} A_{q,\tilde q} := \lim_{t\to\ii} \E\Bigg[ \exp\bigg(\rt1 \bigg[ \sum_{\ell =1}^{\nu_\ii(\Gamma)} q_\ell\,t\1 M^{\ell,r}_t   +  \sum_{k=1}^{2{\bf g}(\Gamma)} \tilde q_k\,t^{-1/2} \tilde M^k_t \bigg]\bigg)  \Bigg]  , \eeq 
which by Item 1) above is the quantity to calculate to get the asymptotic law of Theorem \ref{the.CV}. \  Items 3) and 4) above show that we have : 
\beq \label{f.loilim1} A_{q,\tilde q} = \lim_{r\to\ii}\, \lim_{t\to\ii} \E\Bigg[ \exp\bigg(\rt1 \bigg[ \sum_{\ell =1}^{\nu_\ii(\Gamma)} q_\ell\,t\1 \sum_{n=1}^{[\rr^\ell_r\,t] } \f_n^\ell  +  t^{-1/2} \sum_{n=1}^{[\rr_r\,t]} J_n^q \bigg]\bigg)  \Bigg] . \eeq 
   
   Let us apply now the Markov property : conditionally on the $\,\sigma$-field $\,\FF\,$ generated by the induced Markov chain $(Z_{\sigma_n},Z_{\tau_{n+1}})$, the random variables $\{ \f_n^\ell\,, J_n^q\,|\,1\le \ell\le \nu_\ii(\Gamma) , n\in\N^*\}$ are independent. \  Therefore, denoting by $\,\E^{\FF}$ the conditional expectation with respect to $\,\FF\,$, we have : 
\beq \label{f.loilim2} A_{q,\tilde q} = \lim_{r\to\ii}\, \lim_{t\to\ii} \E\Bigg[  \prod_{\ell =1}^{\nu_\ii(\Gamma)}  \prod_{n=1}^{[\rr^\ell_r\,t] } \E^{\FF}\Big[ e^{\rt1 (q_\ell/t) \f_n^\ell} \Big] \times \prod_{n=1}^{[\rr_r\,t]} \E^{\FF}\bigg[ e^{\rt1 J_n^{q/\sqrt{t}}\,} \bigg] \Bigg] . \eeq 
   
   6) \   We must finally get rid of the conditioning on $\,\FF\,$. \  To do this, we work on each $\, \E^{\FF}\Big[ e^{\rt1\, (q_\ell/t) \f_n^\ell} \Big] $, depending on a single excursion in a given solid cusp ;  to analyse such quantity, we can drop for a while the irrelevant index $\,\ell\,$, and suppose that the width $\,h_\ell\,$ of the cusp is 1, for the sake of notational simplicity.  \  
Now by Section \ref{sec.Brown}, during each excursion near the cusp, we have \  ${\ds \tilde x_s = \tilde x_0 + B\Big( \int_0^s \tilde y_t^2\,dt\Big)}$, for some Brownian motion $(B_s)$ independent from the height component $(\tilde y_s)$. Set \  ${\ds Y:= \int_{\tau_1}^{\sigma_1} y_t^2\,dt\,}$. We have \  $\E\Big[ e^{\rt1 q \,B(Y)}\Big] = \E\Big[ e^{-q^2\, Y/2}\Big] = e^{-|q|\,\sqrt{r}}\,$.  \   
Then for any real $\,q\,$ and any $\,n\in\N^*$ we have : 
$$ \E^{\FF}\Big[ e^{\rt1\, q\, \f_n^\ell} \Big] = \E\Big[ e^{\rt1\, q\,  \f_n^\ell} \Big|\, Z_{\sigma_n},Z_{\tau_{n+1}}\Big] = \E\Big[ e^{\rt1 q \,B(Y)}\Big|\, B(Y)\; {\rm modulo}\; 1 \Big]  \, . $$ 
We have thus to make sure that the knowledge of the value of $\,B(Y)\; modulo\; 1\,$ will perturb the law of $\,B(Y)$ only in a negligible way. \  For this, let us fix $\,u\in\R$ and $\e>0$, and write : 
$$ \E\Big[ e^{\rt1 q\, B (Y) }\,\Big|\, B(Y)\in\, ]u,u+\e[+\Z\Big] - 1 = 
{ \E\bigg[ \Big( e^{\rt1 q\, B(Y)} - 1\Big)\times {\ds\sum_{k\in\sZ}}\, 1_{\{ u < B(Y) -k < u+\e\}}\bigg]
\over \E\bigg[ {\ds\sum_{k\in\sZ}}\, 1_{\{ u < B(Y) -k < u +\e\}}\bigg] } $$ 
$$ =\; { \E\bigg[ (2\pi Y)^{-1/2}{\ds\int_{u}^{u +\e} \sum_{k\in\sZ}} \Big( e^{\rt1 q\, (x+k)} - 1\Big)\,
e^{-(x+k)^2/(2Y)}\, dx\bigg] \over \E\bigg[ (2\pi Y)^{-1/2}{\ds\int_{u}^{u +\e} \sum_{k\in\sZ}} \,
e^{-(x+k)^2/(2Y)}\, dx\bigg] } \; . $$ 
Then observing that \ ${\ds \sup\Big\{ |e^{\rt1 q\,k}-1|\times e^{-k^2/(2Y)}\,\Big|\, k\in\R\Big\} \le   |q|\sqrt{Y}\; }$, we can replace the Riemannian sum above by a Riemannian integral + an error term, in order to get : 
$$ \sum_{k\in\sZ} \Big( e^{\rt1 q\, (x+k)} - 1\Big)\, e^{-(x+k)^2/(2Y)}   = \Big( e^{-\, q^2 Y/2} - 1\Big)\, \sqrt{2\pi Y}\; + \O \Big( |q|\sqrt{Y}\,\Big) \, . $$ 
Hence we have : 
$$ \E\Big[ e^{\rt1 q\, B (Y) }\,\Big|\, B(Y)\in\, ]u,u+\e[+\Z\Big] - 1 =  {\E\Big[ e^{-\, q^2 Y/2} - 1+ \O (|q|) \Big] \over \Big( 1 + \O (\E (Y^{-1/2}))\Big)}  = { e^{-|q|\sqrt{r}} - 1  + \O (|q|)\over  1 + \O (1/\sqrt{r}\,)} \, , $$  
whence : 
\beq \label{f.apprmod} \E^{\FF}\Big[ e^{\rt1\, q\, \f_n^\ell} \Big] = 1- \Big(1 + \O_n (1)/\sqrt{r}\,\Big) |q| \sqrt{r} \, ,  \eeq 
for some uniformly bounded function $\,\O_n(1)\,$ of $\,Z_{\sigma_{n}}\,$.  \par  \smallskip 
   
   7) \   To conclude the proof of Theorem \ref{the.CV}, we note that by 
Birkhoff's ergodic Theorem applied to the Markov chain $\, (Z_{\sigma_{n}})\,$
(via the sequence $\O_n(1)$), Formula (\ref{f.apprmod}) implies : 
$$ \prod_{n=1}^{[\rr^\ell_r\,t] } \E^{\FF}\Big[ e^{\rt1 (q_\ell/t) \f_n^\ell} \Big] = \exp\bigg( - {|q_\ell | \sqrt{r} \over t}\sum_{n=1}^{[\rr^\ell_r\, t]} (1 + \O_n(1)/\sqrt{r}\,) + o(1)\bigg) $$ 
$$ \buildrel{t\to\ii}\over{\lra}\, \exp\Big( - {|q_\ell | \sqrt{r}\,\rr^\ell_r}\, (1 + \O(1/\sqrt{r}\,))\Big)\, . $$ 
Hence we get from Formula (\ref{f.loilim2}) : 
$$ A_{q,\tilde q}  = \lim_{r\to\ii}\,\lim_{t\to\ii}\, \E\Bigg[ \prod_{n=1}^{[\rr_r\,t]} \E^{\FF}\bigg[ e^{\rt1 J_n^{q/\sqrt{t}}\,} \bigg] \Bigg] \times \exp\bigg( -  \sum_{\ell =1}^{\nu_\ii(\Gamma)} {|q_\ell | \sqrt{r}\,\rr^\ell_r}\, (1 + \O(1/\sqrt{r}\,))\bigg) $$ 
$$ = \,\lim_{t\to\ii}\, \E\Bigg[\sum_{k=1}^{2{\bf g}(\Gamma)} \tilde q_k\,t^{-1/2} \tilde M^k_t \bigg]\bigg)  \Bigg] \times \exp\bigg( -  \sum_{\ell =1}^{\nu_\ii(\Gamma)} {|q_\ell | \,h_\ell  \over 2\,V(\Gamma\sm\H^2)} \bigg) $$ 
$$ =\, \exp\bigg( - \sum_{k=1}^{2{\bf g}(\Gamma)} \5\,\tilde q_k\, \langle\tilde\omega_k , \tilde\omega_k\rangle - \sum_{\ell =1}^{\nu_\ii(\Gamma)} {|q_\ell | \,h_\ell  \over 2\,V(\Gamma\sm\H^2)} \bigg) , $$ 
by Item 4), Formula (\ref{f.nbreexc}), and Lemma \ref{lem.ctl}. This achieves the proof of Theorem \ref{the.CV}, owing to Formula (\ref{f.loilim}) defining $\,A_{q,\tilde q}\,$ and to Item 1). $\;\diamond$  

\section{Proof of Theorem \ref{the.geo}} \label{sec.proofg}  \indf 
   The strategy for this proof is mainly to replace the geodesic paths by the Brownian paths, as in [E-LJ], [F2], [LJ1], in order to reduce Theorem \ref{the.geo} to Theorem \ref{the.CV}. As in [F3],  we shall here take advantage of the closedness of the forms $\,\omega_j,\tilde\omega_l\,$, somewhat as in [LJ1], to get a simple enough proof, without using a spectral gap, nor rising to the stable foliation ; simultaneously disintegrating the Liouville and the Wiener measures, we condition the Brownian motion (starting from a given point $z\in\H^2$) to exit the hyperbolic plane at the same point as a given geodesic (starting also from $z$). We essentially follow ([F3], Section 14). \par\smallskip 
   Because of Lemma \ref{lem.w12}, the asymptotic law we are looking for is given by the asymptotic behavior, as $t\to\ii\,$ and for
$\,(\lambda',\lambda)\in\R^{\nu_\ii(\Gamma)}\times\R^{2{\bf g}(\Gamma)}$, of the following quantity :  
\beq \label{f.Jllamb} J_t^{\lambda} := \int_{\Gamma\sm\! G} \exp\left[ \rt1\left( \sum_{j=0}^{\nu_\ii(\Gamma)-1} {\lambda'_j\over t}\int_{g[0,t]}\omega_j'+ \sum_{l=1}^{2{\bf g}(\Gamma)} {\lambda_l\over \sqrt{t}}\int_{g[0,t]}\tilde\omega_l\right)\right] \mu^\Gamma(dg) \eeq
$$ =  {\rm covol}(\Gamma)\1 \int_{\Gamma\sm\! \sH^2} \int_0^{2\pi} \exp\Big[ \rt1\Big(t\1\int_{g(y,x,\t )[0,t]} \omega' + t^{-1/2}\int_{g(y,x,\t)[0,t]} \tilde\omega\Big)\Big]  d\t\,{ dx\, dy\over 2\pi\, y^2}\; , $$ 
where we use the notations of Section \ref{sec.Iwas} and set : 
\beq \label{f.omegas} \omega' := \sum_{j=0}^{\nu_\ii(\Gamma)-1} \lambda'_j\,\omega'_j\; ,\quad  \tilde\omega :=  \sum_{l=1}^{2{\bf g}(\Gamma)} \lambda_l\, \tilde\omega_l\, . \eeq 

\subsection{Conditioning by end-points} \label{sec.endp} \indf
   For $ (z=x+\rt1 y\, ,\, \t )\in\H^2\times (\R /2\pi\Z )\,$, denote by $(z^\t_t)$ the geodesic of $\H^2$ defined by $\,g(z,\t )$, and by $\,\P_z^\t\,$ the law of the Brownian motion $(Z^\t_t)$ of $\H^2$, started from $z$ and conditioned to exit $\H^2$ at the positive end $\,z^\t_\ii\,$ of $(z^\t_t)$. \par \smallskip 

   Consider then the the hitting time, say $\,h_t\,$, by the coordinate process $(Z_t)$, of the stable horocycle defined by $(z^\t_\ii,z^\t_t)$. \  It is defined precisely by \parn 
${\ds h_t = h_t^{z,\t} := \inf\{s>0\,|\, B_{z^\t_\ii}(z,Z_s) = e^t\}}\,$, \quad where \quad 
${\ds (z,z')\mapsto B_u(z,z') = p(z',u)/p(z,u)}\;$ denotes the Busemann function based at 
${\ds u\in\p\H^2\,}$, $\;p$ denoting the Poisson kernel.  \par \smallskip 

   The following lemma ensures that the disintegration of the Liouville and Wiener measures is simultaneous, by conditioning with respect to the end-point $z^\t_\ii\,$. A reason is that the harmonic measures at $\p\H^2$ are the same for both, namely $\, p(z,u)du\,$.  
\begin{lem} \label{lem.Pmu} \quad $\,{\ds \P_z := \int_0^{2\pi}\P_z^\t\,\, {d\t\over 2\pi}}$ \ is the
Wiener measure started from $z$, \  for any \parn $z\in\Gamma\!\sm\! \H^2$, and \quad ${\ds \P_\mu := \int\P_z^\t\, d\mu^\Gamma (z,\t )}$ \   is the stationary Wiener measure on $\,\Gamma\!\sm\! \H^2$. 
\end{lem}
\ub{Proof} \quad  $(Z^\t_t)$ is by definition the $h$-process of the unconditioned Brownian motion, with $\,{\ds h(z)=p(z,z^\t_\ii )}\,$, $\; p(z,u) = y/|z-u|^2\,$ still denoting the Poisson kernel. \parn 
Hence we have for any $\,(z,\t )\,$, any $t$ and any $\FF_t$-measurable positive functional $F_t$ : \parn
\centerline{${\ds \E_z^\t [F_t] = \E_z [B_{z^\t_\ii}(z,Z_t)\times F_t]}\,$.} \parn 
The first identity of the lemma follows, since  for any $z,\t ,Z\,$ we have \parn 
\centerline{${\ds\int_0^{2\pi} B_{z^\t_\ii}(z,Z)d\t = 2\int_{\sR} B_{u}(z,Z)p(z,u)du = 2\int_{\sR} p(Z,u)du = 2\pi}\,$.} \par 
   Integrating this first identity with respect to the normalized Liouville measure $\,\mu^\Gamma$ gives immediately the second identity of the lemma. $\;\diamond $ 

\subsection{From geodesics to Brownian paths} 
\indf 
   We perform here the substitution of the Brownian paths for the geodesics. Our first aim is to establish the following, to the proof of which this section is devoted. 
\begin{prop}\label{pro.contour}\quad  As $t\to\ii\,$, $\,J_t^{\lambda}$ (defined by Formula 
(\ref{f.Jllamb})) behaves as 
$$ K_t^{\lambda',\lambda} := {\rm covol}(\Gamma)\1 \int_{\Gamma\sm\! \sH^2} \int_0^{2\pi} \E_z^\t\left( \exp\Big[ {\ts{\rt1\over t}} \int_{z}^{Z_{h_t}}\omega' + {\ts{\rt1\over\sqrt{t}}}  \int_{z}^{Z_{h_t}} \tilde\omega \Big]\right)  d\t\,{ dx\, dy\over 2\pi\, y^2}\; . $$ 
\end{prop}
\par 
    The forms $\,\omega', \tilde\omega\,$ being closed, we have the following expression for $\,J_t^{\lambda}$ : 
$$ J_t^{\lambda} = \int_{\Gamma\sm\!\sH^2}\! \int_0^{2\pi} \E_z^\t\left( \exp\Big[{\ts{\rt1\over
t}} \Big(\int_{z}^{Z_{h_t}}\! \omega'+\int_{Z_{h_t}}^{z^\t_{t}} \omega'\Big)+ {\ts{\rt1\over\sqrt{t}}}\Big(
\int_{z}^{Z_{h_t}} \tilde\omega +\int_{Z_{h_t}}^{z^\t_{t}}\! \tilde\omega\Big)\Big]\right) 
{ d\t\, dx\, dy\over 2\pi {\rm covol}(\Gamma)  y^2}\, . $$ 

   Applying the isometry $\,f_{z,\t}\,$ of $\H^2$ which maps $\,g(1,0)\,$ to $\,g(z,\t )\,$, we see that 
the law of \  ${\ds \int_{Z_{h_t}}^{z^\t_{t}} \tilde\omega}$ \  under $\,\P_z^\t\,$ is the same as the law of \   ${\ds \int_{Z^0_{h_t}}^{e_t}f_{z,\t}^*\tilde\omega}\,$, where $\,e_t:= \rt1 e^t\,$ and $\,Z^0_{h_t}\,$ is the point at which  the Brownian motion $(Z^0_t)$ started from $\rt1$ and conditioned to exit at $\ii\,$ hits the horizontal horocycle having equation $\,y=e^t\,$. \par \smallskip 
   Now $(Z^0_t)$ is the $h$-process of the unconditioned Brownian motion, with \parn 
$h(z)=p(z,\ii ) \equiv y\,$, so that its infinitesimal generator is $\,\5y\1\D\circ y = \5\D + y\p_y\,$, \parn $\D\,$ denoting the Laplacian of $\H^2$. Thus we have $\;{\ds Z^0_t = \rt1 e^{w_t+t/2} + \int_0^t e^{w_s+s/2}dW_s}\,$, for two independent standard real Brownian motions $(w_t)$ and $(W_t)$. \par 
   As a consequence, using the boundedness of $\,\tilde\omega\,$, we have 
$$ \int_{Z^0_{h_t}}^{e_t}f_{z,\t}^*\tilde\omega = \O\Big( e^{-t}\times\Big|\int_0^{\inf\{s\,|\, w_s+s/2=t\}}e^{w_s+s/2}dW_s \Big|\Big) . $$  

   The technical Brownian behavior we need now and after is given by the following. 
\begin{lem}\label{lem.loi}\quad  As $\,t\to\ii\,$,  $\;{\ds e^{-t}\,\int_0^{\inf\{s\,|\,w_s+s/2=t\}}
e^{w_s+s/2} dW_s }\;$ converges in law, \ and \parn 
$\inf\{s\,|\, w_s+s/2=t\} = 2t + o(t^q)\;$ almost surely, for any $\,q\in ]1/2,1]\,$. 
\end{lem} 
\ub{Proof} \quad Fix $c\in\R\,$, set $\, y^0_t := e^{w_t+t/2}\,$, and look for a $C^2$ function $f$ on $\R_+$ such that \parn 
${\ds R_t:= \exp\Big[{-(c^2/2)\int_0^t (y^0_s)^2 ds}\Big] f(y^0_t)}$ \   be a martingale.
\  $(y^0_t)$ having generator $\,\5 y^2\p_y^2 + y\p_y\,$, we have by It\^o's formula : 
$$ R_t = f(1) + mart + \5\int_0^t e^{-(c^2/2)\int_0^s (y^0_v)^2dv}\times (y^0_s)^2\times\Big[ f''(y^0_s)+
2(y^0_s)\1 f'(y^0_s) -c^2f(y^0_s)\Big] ds\; , $$ 
whence the equation : \quad $ f''(y)+2y\1 f'(y)-c^2f(y) = 0\,$. \ Setting $\; f_1(y):= \sqrt{y} f(y)\,$, 
this gives $\quad f_1''(y)+y\1 f_1'(y)-(c^2+(2y)^{-2})f_1(y) = 0\,$. \ Since $f_1$ must be bounded near $0$,
we have, up to some multiplicative constant : \quad 
${\ds f(y) = (cy)^{-1/2}I_{1/2}(cy) = \sum_{k\ge 0} {(cy)^{2k}\over 2^{2k+\5}k!\G(2k+{3\over 2})}}\,$, where
$\,I_r\,$  denotes the usual modified Bessel function. \par 
   The optional sampling theorem then gives 
$$ \E\Big[\exp\Big(\!\rt1 c \int_0^{\inf\{s\,|\,w_s+s/2=t\}}\! e^{w_s+s/2} dW_s\!\Big)\Big] = 
\E\Big[\exp \Big(\! -{c^2\over 2}\int_0^{\inf\{s\,|\,y^0_s=e^t\}} (y^0_s)^2ds\Big)\Big] = {f(1)\over
f(e^t)} \, . $$ 

   Changing $c$ into $c e^{-t}\,$, we get as $t\to\ii$ : 
$$ \E\Big[\exp (\rt1 c\, e^{-t}\int_0^{\inf\{s\,|\,w_s+s/2=t\}} e^{w_s+s/2} dW_s\Big] \lra 
\Big(\sum_{k\ge 0} {\G (3/2)\, c^{2k}\over 2^{2k}k!\G(2k+{3\over 2})}\Big)\1\; \in L^2 (\R ,dc)\; , $$ 
which proves the first sentence of the lemma. \par 
   Finally, the second sentence of the lemma is straightforward from the following observation : \ 
setting again $\;{\ds h_t = h_t^{\rt1 ,0} = \inf\{s\,|\, w_s+s/2=t\} = \inf\{s\,|\, y^0_s=e^t\}}\;$, we have
\parn
\centerline { ${\ds t = \log\, y^0_{h_t} = \5 h_t + w_{h_t} = \5 h_t + o((h_t)^q)\; }.\;\;\diamond $ }
\medskip 

   As a consequence of this lemma and of the above, we see that \  ${\ds t^{-1/2} \int_{Z_{h_t}}^{z^\t_{t}} \tilde\omega}$ \  goes to 0 in $\,\P_z^\t\,$-probability. \ This proves half of
Proposition \ref{pro.contour}.  \par\smallskip  

   We have now to deal with the law of \  ${\ds t\1\int_{Z_{h_t}}^{z^\t_{t}} \omega'}$ \   under $\,\P_z^\t\,$, or equivalently by the same reason as above for $\,\tilde\omega\,$, with the law of \  ${\ds t\1\int_{Z^0_{h_t}}^{e_t} f_{z,\t}^*\,\omega'}$. This cannot be handled further as above, since $\,\omega'$ is unbounded. But integrating along the horizontal horocycle $\,y=e^t\,$ containing $\,e_t,Z^\t_{h_t}$, we have the following estimate : 
$$ \Big|\int_{Z^0_{h_t}}^{e_t}f_{z,\t}^*\,\omega'\Big| \le \Big| e^{-t}\int_0^{h_t} e^{w_s+s/2}
dW_s\Big| \times \sup\Big\{ |f_{z,\t}^*\,\omega'|_{(\rt1 +x)e^t}\,;\, |x|\le \Big| e^{-t}\int_0^{h_t}
e^{w_s+s/2} dW_s\Big|\Big\} , $$ 
where again \  $ h_t = h_t^{\rt1 ,0} = \inf\{s\,|\, y^0_s=e^t\} = \inf\{s\,|\, w_s+s/2=t\}$. 

   Fix any $r>0\,$. \  Lemma \ref{lem.loi} shows that the laws of $\;{\ds e^{-t}\int_0^{h_t} e^{w_s+s/2} dW_s}\,$, for large $\,t\,$, are tight, and then provides some $R>0\,$ such that 
$\;{\ds\P\Big[\Big| e^{-t}\!\int_0^{h_t} e^{w_s+s/2}dW_s\Big| > R\Big] < r}\,$ for any large enough positive  $t\,$. \par 
 
   We deduce from these last two estimates that 
$$ \P_z^\t\Big[\Big| t\1\int_{Z_{h_t}}^{z^\t_{t}}\omega'\Big| >r\Big] = \P\Big[\Big|
t\1\int_{Z^0_{h_t}}^{e_t}f_{z,\t}^*\,\omega'\Big| >r\Big] \le r + 1_{ \Big\{ t\1\sup\Big\{
|f_{z,\t}^*\,\omega'|_{(\rt1 +x)e^t}\Big| |x|\le R\Big\} > r/R\Big\} } , $$   and then by integrating against
$\mu\,$ and using Lemma \ref{lem.Pmu} : 
$$ \P_\mu\Big[\Big| t\1\int_{Z_{h_t}}^{z^\t_{t}} \omega'\Big| >r\Big] \le r + \mu \Big[ t\1\sup\Big\{
|\omega'|_{H_x(z^\t_t)}\,\Big|\, |x|\le R\Big\} > r/R\Big]  $$
$$ \hskip 39mm = r + \mu \Big[ t\1\sup\Big\{ |\omega'|_{H_x(z)}\,\Big|\, |x|\le R\Big\} > r/R\Big] \, , $$ 
where $\, (H_x\, ,\, x\in\R)\,$ denotes the positive horocycle flow. For the last equality, we used the
invariance of the Liouville measure $\mu$ under the geodesic flow. \parn 
 By continuity of $\,|\omega'|\,$, 
$\,\sup\Big\{ |\omega'|_{H_x(z)}\Big| |x|\le R\Big\}\,$ is finite for every $z$, and thus we just proved :    
$$ \P_\mu\Big[\Big| t\1\int_{Z_{h_t}}^{z^\t_{t}} \omega'\Big| >r\Big] \le 2r \quad \hbox{ for large enough }
t \, . $$ 

   Since in the last expression above for $\,J_t^{\lambda}$ (immediately after Proposition \ref{pro.contour}), we were not only under the law $\,\P_z^\t\,$, but indeed under the law $\P_\mu = \bigint \P_z^\t\, d\mu (z,\t )\,$, 
we have so far proved Proposition \ref{pro.contour}. 

\subsection{End of the proof of Theorem \ref{the.geo}} 
\indf 
   Section \ref{sec.Brown} allows to denote also by $\P_\mu$ the stationary Wiener measure on $\,\Gamma\!\sm\! G\,$, since the Brownian motion of $\,G\,$ projects on the Brownian motion of $\,\H^2$ (and similarly for the volume measures).  Recall also that the forms $\,\omega'_j,\tilde\omega_l\,$ come from $\,\Gamma\!\sm\!\H^2$ : they are defined on $\,\Gamma\!\sm\! G\,$ and on $\,\Gamma\!\sm\!\H^2$ as well, in other words are invariant under pull back $\pi^*$ by the canonical projection. Hence the joint laws of their integrals along the Brownian paths are the same, no matter whether they are understood on $\,\Gamma\!\sm\! G\,$ or on $\,\Gamma\!\sm\!\H^2$. \par 
   Moreover we have seen in Section \ref{sec.Brown} also that the angular Brownian component $\,\t_s\,$ is a mere one-dimensional Brownian motion. As a consequence, it is immediate that \quad 
${\ds t\1\int_{g[0,t]}d\t = (\t_t-\t_0)/t}\;$ goes to 0 $\,\P_\mu$-almost surely. \ Therefore we can replace in Theorem \ref{the.CV} the form $\o_0\,$ by the form $\,\o'_0 = \o_0 -d\t\,$. \par 
  These remarks show that the following is merely an alternative version of Theorem \ref{the.CV} 
(with the notations of Formula (\ref{f.omegas}) and Theorem \ref{the.CV}). 
\begin{cor} \label{cor.CV} \quad We have for any $\,(\lambda',\lambda)\in\R^{\nu_\ii(\Gamma)}\times\R^{2{\bf g}(\Gamma)}$ : 
$$ \lim_{t\to\ii}\, \E_\mu\left( \exp\Big[ {\ts{\rt1\over t}} \int_{Z[0,t]} \omega'
+ {\ts{\rt1\over\sqrt{t}}} \int_{Z[0,t]} \tilde\omega \Big]\right) $$
$$  = \Lambda(\lambda',\lambda) := \E\left( \exp\left[ \rt1 \left( \sum_{j=0}^{\nu_\ii(\Gamma)-1} \lambda'_j  \sum_{\ell=1}^{\nu_\ii(\Gamma)} r^\ell_j\, {\mathcal Q}_\ell + \sum_{l=1}^{2{\bf g}(\Gamma)} \lambda_l\, {\mathcal N}^l \right)\right]\right) . $$ 
\end{cor}

   Now Lemma \ref{lem.loi} asserts that the time-change $\,h_t= h_t^{z,\t}$ appearing in the expression of $\,K_t^{\lambda',\lambda}$ in Proposition \ref{pro.contour}, satisfies \   $h_t = 2t + o(t)$ $\,\P_z^\t$-almost surely, uniformly with respect to $(z,\t )\,$. Indeed, the law under $\,\P_z^\t\,$ of $\,h_t^{z,\t}$ equals the law of $\,h_t^{\rt1 ,0}\,$ in Lemma \ref{lem.loi}. \   So that, with arbitrary large probability, we can write \  $h_t = 2t + o(t)$ \   with a uniform deterministic $o(t)\,$. \      This allows to replace $\,t\,$ by $\,h_{t}\,$ in the formula of Corollary \ref{cor.CV} above, getting then (using also the definition of $\,\P_\mu\,$ in Lemma \ref{lem.Pmu}) : 
$$ \lim_{t\to\ii}\, K_{t}^{\lambda'/2,\lambda/\sqrt{2}} = \lim_{t\to\ii}\, \E_\mu\left( \exp\Big[ {\ts{\rt1\over 2t}} \int_{Z[0,h_t]} \omega'
+ {\ts{\rt1\over\sqrt{2t}}} \int_{Z[0,h_t]} \tilde\omega \Big]\right) = \Lambda(\lambda',\lambda) . $$ 
Therefore, using Proposition \ref{pro.contour} \  we have proved that 
$$ \lim_{t\to\ii} J_t^{\lambda} = \Lambda(2\lambda',\sqrt{2}\,\lambda) . $$  

   This concludes the proof, since by the definition of $\, \Lambda\,$ in Corollary \ref{cor.CV}, by the very definition (\ref{f.Jllamb}) of $J_t^{\lambda}$, and by Lemma \ref{lem.w12}, this formula is equivalent to Theorem \ref{the.geo}. 

\bigskip \bigskip \bigskip 
\centerline{{\Large{\bf REFERENCES}}}  
\bigskip \bigskip 

\vbox{ \noindent 
{\bf [D]} \  DieudonnŽ J. \  \  {\it \'ElŽments d'Analyse 9.} \hskip 17mm   Gauthier-Villars, Paris, 1982.} 
\bigskip 

\vbox{ \noindent 
{\bf [E-F-LJ1]} \ Enriquez N. , Franchi J. , Le Jan Y.  \ \  {\it Stable windings on hyperbolic surfaces.
} \par \smallskip \hskip 25mm    Prob. Th. Rel. Fields 119, 213-255, 2001. }
\bigskip 

\vbox{ \noindent 
{\bf [E-F-LJ2]} \ Enriquez N. , Franchi J. , Le Jan Y.    \ 
{\it Central limit theorem for the geodesic \par \hskip 16mm flow associated with a Kleinian group, 
case  $\delta > d/2$.} \par 
\hskip 25mm   J. Math. Pures Appl. 80, 2, 153-175, 2001. } 
\bigskip 

\vbox{ \noindent 
{\bf [E-LJ]} \  Enriquez N. , Le Jan Y.  \quad {\it Statistic of the winding of geodesics on 
a Riemann \par \hskip 10mm surface with finite volume and constant negative curvature.}  \par 
\hskip 25mm {Revista Mat. Iberoam.}, vol. 13, n$^o\, 2$, 377-401, 1997. }
\bigskip 

\vbox{ \noindent 
{\bf [E-MM]} \  Eskin A. , McMullen C.  \quad {\it Mixing, counting, and equidistribution in Lie groups.}  \par 
\hskip 25mm {Duke Math. J.}, vol. 71, n$^o\, 1$, 181-209, 1993. }
\bigskip 

\vbox{ \noindent 
{\bf [F1]} \  Franchi J. \ {\it Asymptotic singular windings of ergodic diffusions.}
\smallskip \par 
\hskip 26mm   Stoch. Proc. and their Appl., vol. 62, 277-298, 1996. } 
\bigskip 

\vbox{ \noindent 
{\bf [F2]} \  Franchi J. \ {\it Asymptotic singular homology of a complete hyperbolic 3-manifold 
of finite \par \hskip 4mm volume.} 
\smallskip 
\hskip 5mm  Proc. London Math. Soc. (3) 79, 451-480, 1999.  } 
\bigskip 

\vbox{  \noindent 
{\bf [F3]} \  Franchi J. \  {\it Asymptotic windings over the trefoil knot.} 
\par \hskip 25mm  Revista Mat. Iberoam., vol. 21, n$^o\, 3$, 729-770, 2005.  } 
\bigskip 

\vbox{  \noindent 
{\bf [G-LJ]} \  Guivarc'h Y. , Le Jan Y. \  \ {\it  Asymptotic windings of the geodesic flow on modular \par \hskip 5mm    surfaces with continuous fractions.} \    Ann. Sci. \'Ec. Norm.
Sup. 26, n$^o$4, 23-50, 1993.   }
\bigskip 

\vbox{ \noindent 
{\bf [H]} \ Hopf E. \ \ {\it Ergodicity theory and the geodesic flow on a surface of constant negative \par  \hskip 21mm    curvature.} \hskip 12mm   Bull. Amer. Math. Soc. 77, 863-877, 1971.} 
\bigskip 

\vbox{  \noindent 
{\bf [I-W]} \ Ikeda N. ,  Watanabe S. \quad {\it Stochastic differential equations and diffusion processes.}
\par \hskip 27mm North-Holland Kodansha, 1981. }
\bigskip 

\vbox{ \noindent 
{\bf [L]} \  Lehner  J.  \quad  {\it Discontinuous groups and automorphic functions.}  \par 
\hskip 25mm    Amer. Math. Soc, math. surveys n$^o$VIII, Providence, 1964. } 
\bigskip 

\vbox{ \noindent 
{\bf [LJ1]} \  Le Jan  Y.  \quad  {\it Sur l'enroulement g\'eod\'esique des surfaces de Riemann.}  \par 
\hskip 30mm    C.R.A.S. Paris, vol 314, S\'erie I, 763-765, 1992. } 
\bigskip 

\vbox{ \noindent 
{\bf [LJ2]} \  Le Jan  Y.  \quad  {\it The central limit theorem for the geodesic flow on non compact 
 \par\hskip 8mm manifolds of constant negative curvature. }  \hskip 1mm    Duke Math. J.  (1) 74, 159-175, 1994.  } 
\bigskip 

\vbox{ \noindent 
{\bf [M]} \ Miyake T. \ \ {\it Modular forms.} \quad  Springer, Berlin 1989.} 
\bigskip 

\vbox{ \noindent 
{\bf [R]} \ Randol B.   \ {\it The behavior under rojection of dilating sets in a covering space.} \par    \hskip 25mm  Trans. Amer. Math. Soc. 285, 855-859, 1984. }
\bigskip 

\vbox{ \noindent 
{\bf [R-Y]} Revuz D. , Yor M.   \ {\it Continuous martingales and Brownian motion.} \   Springer,  1999. }
\bigskip 

\if{ \vbox{  \noindent 
{\bf [SL-M]} de Sam Lazaro J. , Meyer P.A. \quad {\it Questions de th\'eorie des flots.}  \par 
\hskip 23mm  S\'em. Probab. IX, Lect. Notes n$^o\,$465, P.A. Meyer editor, Springer 1975. }
\bigskip } \fi 

\vbox{ \noindent 
{\bf [Sh]} \ Shimura G. \quad {\it Introduction to the arithmetic theory of automorphic functions.} \par 
\hskip 31mm   Publ. Math. Soc. Japan,  Princeton University Press, 1971. }
\bigskip 

\vbox{ \noindent 
{\bf [Sp]} \ Spitzer F. \quad {\it Some theorems concerning two-dimensional Brownian motion.} \par 
\hskip 28mm    Trans. A. M. S. vol. 87, 187-197, 1958. }
\bigskip 

\vbox{ \noindent 
{\bf [T]} \  Thurston W.P. \ {\it Three dimensional manifolds, Kleinian groups and  
hyperbolic \par \hskip 4mm geometry.}  \quad { Bull. Amer. Math. Soc.} 6, 357-381, 1982. }
\bigskip 

\vbox{ \noindent 
{\bf [W]} \  Watanabe S. \ {\it Asymptotic windings of Brownian motion paths  on Riemannian 
\par \hskip 4mm surfaces. }  \quad {\, Acta Appl. Math.} 63, n$^o\,$1-3, 441-464, 2000. }
\bigskip

\end{document}